\documentclass[10pt]{article}
\usepackage{latexsym}
\usepackage{amsbsy}
\usepackage{amssymb}
\usepackage{amsmath}
\usepackage[latin1]{inputenc}
\usepackage{a4wide}
 \newtheorem{thm}{Theorem}[section]
 \newtheorem{defin}[thm]{Definition}
 \newtheorem{lem}[thm]{Lemma}
 \newtheorem{prop}[thm]{Proposition}
 \newtheorem{cor}[thm]{Corollary}
 \newtheorem{rem}[thm]{Remark}
 \newtheorem{ex}[thm]{Example}

 \newcommand{\bthm}{\begin{thm}}
 \newcommand{\ethm}{\end{thm}}
 \newcommand{\bd}{\begin{defin}}
 \newcommand{\ed}{\end{defin}}
 \newcommand{\blem}{\begin{lem}}
 \newcommand{\elem}{\end{lem}}
 \newcommand{\bcor}{\begin{cor}}
 \newcommand{\ecor}{\end{cor}}
 \newcommand{\bprop}{\begin{prop}}
 \newcommand{\eprop}{\end{prop}}
 \newcommand{\brem}{\begin{rem} \rm}
 \newcommand{\erem}{\end{rem}}
 \newcommand{\bex}{\begin{ex} \rm}
 \newcommand{\eex}{\end{ex}}

 \newcommand{\pr}{\noindent{\bf Proof. }}
 \newcommand{\ep}{\nolinebreak{\hspace*{\fill}$\Box$ \vspace*{0.25cm}}}

 \newcommand{\beq}{ \begin{equation} }
 \newcommand{\eeq}{\end{equation} }

 \newcommand{\bea}{\begin{eqnarray}}
 \newcommand{\eea}{\end{eqnarray}}

 \newcommand{\beas}{\begin{eqnarray*}}
 \newcommand{\eeas}{\end{eqnarray*}}

 \newcommand{\beqs}{\begin{equation*}}
 \newcommand{\eeqs}{\end{equation*}}

 \newcommand{\bi}{\begin{itemize}}
 \newcommand{\ei}{\end{itemize}}

 \newcommand{\ben}{\begin{enumerate}}
 \newcommand{\een}{\end{enumerate}}

 \newcommand{\ba}{\begin{array}}
 \newcommand{\ea}{\end{array}}

 \newcommand{\ds}{\displaystyle}

 \newcommand{\R}{\mathbb R}
 \newcommand{\N}{\mathbb N}
 \newcommand{\C}{\mathbb C}

 \newcommand{\cA}{\ensuremath{{\mathcal A}}}
 \newcommand{\cC}{\ensuremath{{\mathcal C}}}
 \newcommand{\cD}{\ensuremath{{\mathcal D}}}
 \newcommand{\cE}{\ensuremath{{\mathcal E}}}
 \newcommand{\cG}{\ensuremath{{\mathcal G}}}

 \newcommand{\cL}{\ensuremath{{\mathcal L}}}
 \newcommand{\cN}{\ensuremath{{\mathcal N}}}

 \newcommand{\gs}{\ensuremath{\mathcal G}}
 
 \newcommand{\pd}{\partial}

 \newcommand{\eps}{\varepsilon}
 \newcommand{\vphi}{\varphi}

 \newcommand{\id}{\ensuremath{\mathrm{id}}}
 \newcommand{\Div}{\ensuremath{\mathrm{Div\,}}}

 \newcommand{\sign}{\ensuremath{\mbox{sign}}}
 \newcommand{\supp}{\ensuremath{\mathrm{supp}}}

 \newcommand{\vf}{{\bf v}}
 
 \newcommand{\un}{u^{(n)}}
 \newcommand{\prn}{{\rm pr}^{(n)}}
 
 \newcommand{\prfirst}{{\rm pr}^{(1)}}

 \newcommand{\tils}{\tilde{s}}
 \newcommand{\tilx}{\tilde{x}}
 \newcommand{\tily}{\tilde{y}}
 \newcommand{\tilR}{\tilde{\R}}
 
 \newcommand{\tilcG}{\tilde{\ensuremath{{\cal G}}}}

\newcommand{\DIV}{\operatorname{\rm{div}}}
\newcommand{\p}{\partial}
\newcommand{\beps}{\mbox{\boldmath$\varepsilon$}}
\newcommand{\bsig}{\mbox{\boldmath$\sigma$}}

\begin{document}

 \title{Foundations of the calculus of variations in generalized function algebras}
 \author{Sanja Konjik
         \footnote{Faculty of Agriculture, University of Novi Sad,
         Trg Dositeja Obradovi\' ca 8, 21000 Novi Sad, Serbia,
         Electronic mail: kinjoki@uns.ns.ac.yu}\\
         Michael Kunzinger
         \footnote{Faculty of Mathematics, University of Vienna,
         Nordbergstr.\ 15, A-1090 Wien, Austria,
         Electronic mail: michael.kunzinger@univie.ac.at. Work supported by FWF-project
         P-16742 and START-project Y-237.}\\
         Michael Oberguggenberger
         \footnote{ Arbeitsbereich f\"ur Technische Mathematik,
         Institut f\"ur Grundlagen der Bauingenieurwissenschaften , University of
         Innsbruck, Technikerstr.\ 13, A-6020 Innsbruck, Austria,
         Electronic mail: michael.oberguggenberger@uibk.ac.at}\\
       }

 \date{}
 \maketitle

\begin{abstract}
We propose the use of algebras of generalized functions for the analysis of certain
highly singular problems in the calculus of variations. 
After a general study of extremal problems on open subsets of Euclidean space in this setting we 
introduce the first and second variation of a variational problem. We then derive necessary
(Euler-Lagrange equations) and sufficient conditions for extremals. The concept of association
is used to obtain connections to a distributional description of singular variational problems.
We study variational symmetries and derive an appropriate version of N\"other's
theorem. Finally, a number of applications to geometry, mechanics, elastostatics and elastodynamics
are presented.

\end{abstract}

\section{Introduction}
The study of singular problems in the calculus of variations has a long history and
a sizeable literature on diverse aspects of this topic is available (cf., e.g., 
\cite{graves,davie,tuckey} and the literature cited therein).
In this paper we introduce an approach to variational problems involving singularities
based on the theory of algebras of generalized functions in the sense of Colombeau (\cite{c1,c2}).
We are interested in extremizing functionals which are either singular (e.g., distributional)
themselves or whose set of admissible variations includes generalized functions.

When considering singular variational problems in the above sense, 
distribution theory is only of limited use due to the nonlinear nature
of the typical variational problems. 
The theory of algebras of generalized functions, on the other hand, 
provides a nonlinear extension of distribution theory which allows
to model nonlinear singular problems while at the same time
providing optimal consistency properties with respect to the
linear theory. It has found an increasing number of applications
in linear and nonlinear partial differential equations (cf., e.g., 
\cite{MObook}, or \cite{NPsurvey} for a recent survey), 
regularity theory and microlocal analysis (e.g., \cite{garetto_hoermann, 
HOP, delcroix}), as well as in non-smooth differential geometry 
(e.g., \cite{book,ndg,gprg,connections}). As a rather novel development,
the order structure in Colombeau type algebras of generalized functions
has been investigated in \cite{AJ,OPS-positivity,EbM-AF}. It allows the
formulation of variational problems in the generalized setting and
will therefore be crucial to this work.

Furthermore, we study symmetry properties of variational problems,
thereby continuing our investigations on Lie symmetries of differential equations
in generalized functions (\cite{symm, DKP, KK-ggags, KK-gf04}). 
Our aim here is to derive infinitesimal criteria and establish an
appropriate version of N\"other's theorem.

The paper is organized as follows. After recalling some standard 
notations from symmetry group analysis and algebras of generalized
functions, in Section \ref{diffcalc} we examine extremal values
of Co\-lom\-beau generalized functions on open subsets of Euclidean 
space. As in the classical calculus of variations this will provide
the basis for deriving necessary and sufficient conditions for 
general variational problems later on. Already at this stage, new
effects occur due to the structural properties of the set of scalars 
in Colombeau theory: they form an ordered ring but not a field.
The presence of zero divisors requires a more refined analysis than
in the smooth case. In Section \ref{firstsec} we turn to calculus of variations
proper. After singling out the appropriate spaces of generalized functions for
formulating the problem, we study necessary and sufficient conditions for
minimizers. In particular, we prove a fundamental lemma of the calculus of
variations and obtain the Euler-Lagrange equations as necessary conditions.
We also introduce the concept of minimizer in the sense of associations which
allows to study variational problems on the distributional level. 
In section \ref{varsymmnoether} we investigate variational symmetries
in the framework of algebras of generalized functions. Here we derive an
infinitesimal criterion and establish a N\"other theorem, allowing
to derive conservation laws from variational symmetries.
Finally, section \ref{sec:Appl} is devoted to a variety of applications,
focussing on variational problems with singularities that require formulation
in the Colombeau setting. The first example treats geodesics in a generalized
Riemannian metric. Then we turn to particle mechanics with singular potentials.
The third type of example is concerned with elasticity theory, either with degenerate
material properties (loss of ellipticity) or with singular potentials: membranes with springs attached,
beams with discontinuous or vanishing coefficients and rods with generalized stress-strain
relation. In the last example we turn to the hyperbolic case: the wave equation with singular and
nonlinear potential.

 We now briefly introduce the notations that will be used
 throughout this paper. Our basic reference for the theory of
 algebras of generalized functions is \cite{book}. For symmetry group
 analysis we refer to \cite{Olv}. Concerning the calculus of variations,
 any standard text on the subject will cover what is needed in our
 approach (e.g., \cite{GelfandF,GiaquintaH,Jost})
 
 The spaces of independent and dependent variables will be $\R^p$
 and $\R^q$ respectively, or their subsets $\Omega$ and $U$
 (we set $M=\Omega\times U$). By $(x,\un)$ we
 shall denote a point in the $n$-th order prolongation or jet space
 $M^{(n)}$, whose components are $x^i$,
 $i=1,\dots,p$, of $x$ and $u^\alpha_J$, $\alpha=1,\dots,q$,
 $J=(j_1,\dots,j_k), 1\leq j_k \leq p, 1\leq k\leq n$, representing $u$ and
 all partial derivatives of $u$ up to order $n$. The $n$-th
 prolongation of a function $f:\Omega\to U$ will be the function
 $\prn f:\Omega\to U^{(n)}$ with components representing $f$ and
 its partial derivatives up to order $n$. In general, the $n$-th
 prolongation of any object will be denoted by $\prn$. By $G$ we
 shall denote a one-parameter Lie group of transformations which
 acts on $M$ and write for the action $(x',u')=g_\eta\cdot (x,u) =
 \Phi(\eta,(x,u))= (\Xi(\eta,(x,u)), \Psi(\eta(x,u)))$. To
 simplify notation we will also write $\Phi_\eta(x,u)$,
 $\Xi_\eta(x,u)$ and $\Psi_\eta(x,u)$. $G$ is called projectable
 if $\Xi_\eta(x,u)=\Xi_\eta(x)$. The infinitesimal generator of
 $G$ is a vector field on $M$ denoted by $\vf$
 and written in the form  $\sum_{i=1}^p \xi^i (x,u) \pd_{x^i} +
 \sum_{\alpha=1}^q \phi_\alpha (x,u) \pd_{u^\alpha}$.
 The total derivative with respect to $x^i$ will be denoted by
 $D_i$.

 The algebras of generalized functions on which our approach is based will
 be those of special (or simplified) type (\cite[Ch.\ 1]{book}), so we 
 will henceforth drop the adjective ``special'' and simply refer to them
 as Colombeau algebras. The algebra $\cG(\Omega)$ is the quotient $\cE_M(\Omega)/\cN(\Omega)$,
 where
 \beas
 \cE_M(\Omega) &\!\!\!\!=\!\!\!\!& \{(u_\eps)_\eps\in(\cC^\infty(\Omega))^{(0,1]}|
 \forall K\subset\subset\Omega\ \forall \alpha\in\N^p_0\ \exists
 N\in\N: \sup_{x\in K} |\pd^\alpha u_\eps(x)|=O(\eps^{-N})
 \mbox{ as } \eps\to 0\} \\
 \cN(\Omega) &\!\!\!\!=\!\!\!\!& \{(u_\eps)_\eps\in\cC^\infty(\Omega)^{(0,1]}|
 \forall K\subset\subset\Omega\ \forall \alpha\in\N^p_0 \forall m\in\N:
 \sup_{x\in K} |\pd^\alpha u_\eps(x)|=O(\eps^{m})
 \mbox{ as } \eps\to 0\}.
 \eeas
 $\cN(\Omega)$ is an ideal in the algebra $\cE_M(\Omega)$ (all operations
 are defined componentwise, i.e., for fixed $\eps$).
 Elements of $\cE_M(\Omega)$ and $\cN(\Omega)$ are called
 moderate, resp.\ negligible nets of smooth functions.
 $\cG(\Omega)$ is an associative, commutative, differential
 algebra whose elements are equivalence classes
 $u=[(u_\eps)_\eps]$. $\Omega\to\cG(\Omega)$ is a fine sheaf of
 differential algebras on $\R^p$. In particular, there is a well-defined notion
 of support in $\cG(\Omega)$ and we denote by $\cG_c(\Omega)$ the space of
 compactly supported generalized functions in $\cG(\Omega)$.
 Moreover, $\cG(\Omega)$ contains the space $\cD'(\Omega)$
 as a linear subspace and $\cC^\infty(\Omega)$ is a faithful subalgebra: 
 embedding is effected basically by convolution with a fixed mollifier. This embedding is
 a sheaf morphism that commutes with partial differentiation. The concept of association
 provides a means of assigning macroscopic properties to elements of Colombeau
 algebras: elements $u$, $v$ of a Colombeau algebra are called associated if
 $u_\eps - v_\eps \to 0$ in $\cD'$. $u$ is associated to the distribution $w$
 if $u_\eps \to w$ weakly. These notions are independent of the chosen representatives.
  
 Given another open set $\Omega' \subseteq \R^q$, we may consider those elements
 of $\cG(\Omega)^q$ possessing a representative $(u_\eps)_\eps$ such that
 $u_\eps(\Omega)\subseteq \Omega'$ for all $\eps$ and which is compactly bounded
 or c-bounded in the sense that
  $$
 \forall K\subset\subset \Omega\ \exists K'\subset\subset \Omega'\
 \exists \eps_0\in(0,1]\ \forall \eps<\eps_0: u_\eps(K)\subseteq
 K'\,.
 $$
 The space of all $c$-bounded generalized functions from
 $\Omega$ to $\Omega'$ is denoted by $\cG[\Omega,\Omega']$. Elements
 of $\cG[\Omega,\Omega']$ can be composed unrestrictedly as well as
 inserted into elements of $\cG(\Omega')$. A similar definition can be
 given for smooth manifolds instead of open sets and a functorial theory
 of generalized functions based on this notion has been developed in \cite{gfvm,gfvm2}.
 
 The algebra of tempered generalized functions $\cG_\tau(\Omega)$ is
 defined as the quotient
 $\cE_\tau(\Omega)/\cN_\tau(\Omega)$, where
 \beas
 \cE_\tau(\Omega) &\!\!\!\!=\!\!\!\!& \{(u_\eps)_\eps\in \cC^\infty(\Omega)^{(0,1]}|
 \forall \alpha\in\N^p_0\ \exists N\in\N: \sup_{x\in \Omega}
 (1+|x|)^{-N}|\pd^\alpha u_\eps(x)|=O(\eps^{-N})
 \mbox{ as } \eps\to 0\} \\
 \cN_\tau(\Omega) &\!\!\!\!=\!\!\!\!& \{(u_\eps)_\eps\in\cC^\infty(\Omega)^{(0,1]}|
 \forall \alpha\in \N_0^p \ \exists l\in\N\ \forall m\in\N:
 \sup_{x\in \Omega} (1+|x|)^{-l} |u_\eps(x)|=O(\eps^{m})
 \mbox{ as } \eps\to 0\}
 \eeas
 The algebra $\cG_\tau$ can be used to implement Fourier transformation in the Colombeau
 setting (analogous to the space ${\mathcal S}'$ in distribution theory, which in turn is
 embedded into $\cG_\tau$). Our main interest
 in it, however, relies on the fact that insertion of elements of $\cG$ into elements of
 $\cG_\tau$ yields well-defined elements of $\cG$.
 
 We shall also make use of the following 'mixed' variant of Colombeau algebra: let 
 $\tilcG_\tau(\Omega\times\Omega')=
 \tilde{\cE}_\tau(\Omega\times\Omega') /
 \tilde{\cN}_\tau(\Omega\times\Omega')$, where $\Omega'$ is a
 subset of $\R^{p'}$, and
 \beas
 \tilde{\cE}_\tau(\Omega\times\Omega') &\!\!\!\!=\!\!\!\!&
 \{(u_\eps)_\eps\in\C^\infty(\Omega\times\Omega')^{(0,1]}|
 \forall K\subset\subset\Omega\ \forall \alpha\in\N^{p+p'}_0\ \exists
 N\in\N:\\
 && \qquad
 \sup_{x\in K, y\in\Omega'} (1+|y|)^{-N}||\pd^\alpha u_\eps(x,y)|
 =O(\eps^{-N})
 \mbox{ as } \eps\to 0\} \\
 \tilde{\cN}_\tau(\Omega\times\Omega') &\!\!\!\!=\!\!\!\!&
 \{(u_\eps)_\eps\in\C^\infty(\Omega\times\Omega')^{(0,1]}|
 \forall K\subset\subset\Omega\ \forall \alpha\in\N^{p+p'}_0\
 \exists l\in\N\ \forall m\in\N:\\
 && \qquad
 \sup_{x\in K, y\in\Omega'} (1+|y|)^{-l}||\pd^\alpha u_\eps(x,y)|
 =O(\eps^{m})
 \mbox{ as } \eps\to 0\}\,.
 \eeas
 Thus, its elements satisfy $\cG$-estimates in the $\Omega$-variable and
 $\cG_\tau$-estimates in the $\Omega'$-variable.

 It may be seen as a nonstandard feature of Colombeau algebras
 that, due to the presence of infinitesimals, their elements are not uniquely
 determined by their values on the usual points in their respective domains.
 In order to achieve such a unique determination one has to use so-called
 generalized points in the following sense (\cite{point}):
 The space of generalized points $\tilde{\Omega}$ is the set of all equivalence classes in
 $$
 \Omega_M=\{(x_\eps)_\eps\in\Omega^{(0,1]}|\exists N\in\N: |x_\eps|=
 O(\eps^{-N}) \mbox{ as } \eps\to 0 \},
 $$
 where the equivalence relation is defined by
 $$
 (x_\eps)\sim (y_\eps)_\eps \Leftrightarrow \forall m\in\N: |x_\eps
 -y_\eps| =O(\eps^m) \mbox{ as } \eps\to 0,
 $$
 i.e. $\tilde{\Omega}=\Omega_M/\sim$. The set of compactly supported
 points is
 $$
 \tilde{\Omega}_c = \{\tilx=[(x_\eps)_\eps]\in\tilde{\Omega}|
 \exists K\subset\subset \Omega\ \exists\eps_0\in(0,1]\ \forall
 \eps<\eps_0: x_\eps\in K\}.
 $$
 If $\Omega=\R^p$ we shall write $\tilR^p$ resp.\ $\tilR^p_c$. As a
 special case, when $\Omega=\R$, we obtain $\tilR$, the ring of
 constants in any of the above algebras. 
 Elements of $\cG(\Omega)$ are uniquely determined by their point
 values in all compactly supported generalized points in $\Omega$.
 If $\Omega$ is an $n$-dimensional
 box and $u\in\cG_\tau(\Omega)$ then $u=0$ in $\cG_\tau(\Omega)$
 if and only if $u(\tilx)=0$ in $\tilR$, for all $\tilx\in
 \tilde{\Omega}$. Moreover, one can show that Colombeau generalized functions
 are in fact uniquely determined by their values in all near-standard points. Here, a point
 $\tilx\in\tilde{\Omega}_c$ is called near-standard if there
 exists $x\in\Omega$ such that $x_\eps\to x$ as $\eps\to 0$ for
 any representative $(x_\eps)_\eps$ of $\tilx$.

 The above notions can be employed to define generalized Lie group actions
 in the Colombeau setting (cf.\ \cite{symm,DKP,KK-ggags}):
 A generalized group action on $\R^p$ is an element
 $\Phi\in\tilde{\cG}_\tau(\R\times\R^p)^p$ which is a one-parameter group 
 in the sense that $\Phi(0,\,.\,)=\id$ in $\cG_\tau(\R^p)^p$ and
 $\Phi(\eta_1+\eta_2,\cdot)= \Phi(\eta_1,\Phi(\eta_2,\,.\,))$ in
 $\tilde{\cG}_\tau(\R^2\times\R^p)^p$. If $\xi\in\cG_\tau(\R^p)^p$
 is a generalized vector field with the property that there is a
 unique generalized group action
 $\Phi\in\tilde{\cG}_\tau(\R\times\R^p)^p$ satisfying
 $\frac{d}{d\eta}\Phi(\eta,x)=\xi(\Phi(\eta,x))$ in $\tilde\cG_\tau(\R^{1+p})^p$
 (i.e., if $\Phi$ is the unique flow of $\xi$)
 then $\xi$ is called the infinitesimal generator of $\Phi$ and
 both $\xi$ and $\Phi$ are called $\cG$-complete. If, in addition,
 $\Phi$ and $\xi$ have representatives $(\Phi_\eps)_\eps$ and
 $(\xi_\eps)_\eps$ such that $\Phi_\eps$ is the flow of $\xi_\eps$
 for each $\eps\in(0,1]$, then $\Phi$ and $\xi$ are called strictly
 $\cG$-complete.

 The remainder of this introduction will be devoted to the notion of
 positivity in the generalized functions setting. Since the results in this area
 are rather recent and come from a variety of different sources,
 we collect below those parts of the theory which will be crucial for our 
 further investigation (cf.\ \cite{book,OPS-positivity,EbM-AF}). 
 
By \cite[Theorem 1.2.38]{book}, an element $\tilde{x}$ of $\tilde{\R}$ is invertible
if and only if it is {\em strictly nonzero} in the following sense:
for any representative $(x_\eps)_\eps$ of $\tilx$ there exist $a>0$ and $\eps_0$ such that
$|x_\eps|\geq\eps^a$, for all $\eps<\eps_0$. This, in turn, is equivalent to 
$\tilde x$ not being a zero divisor (\cite[Theorem 1.2.39]{book}).
 
 $\tilde{x}\in\tilde{\R}$
 is said to be nonnegative, $\tilde{x}\geq 0$, if there exists a 
 representative such that each $x_\eps$ is nonnegative. Equivalently, for each representative
 $(x_\eps)_\eps$ of $\tilde{x}$ and each $a>0$ there exists
 $\eps_0$ such that $x_\eps +\eps^a\geq 0$, for all $\eps<\eps_0$.
 $\tilde{x}\in\tilde{\R}$ is said to be strictly positive, $\tilde{x}> 0$,
 if it is positive and invertible. By the above, this means that
 there exists $a>0$ such that for each representative $(x_\eps)_\eps$ of
 $\tilde{x}$ there exists $\eps_0$ with $x_\eps >\eps^{a}$, for all $\eps<\eps_0$.
 Note that since $\tilde \R$ contains zero divisors, $\tilde x > 0$ is strictly
 stronger than $\tilde x \ge 0$ and $\tilde x\not=0$.
 We shall need the following characterizations from
 \cite{EbM-AF}: $\tilde{x}\in\tilR$ is strictly positive if and only
 if
 
\begin{itemize}
\item[(i)] $\tilde{x}$ is invertible and has a representative
 $(x_\eps)_\eps$ such that $x_\eps>0, \forall \eps\in(0,1]$;
\item[(ii)] for each representative $(x_\eps)_\eps$ of $\tilde{x}$ there
 exists an $\eps_0$ such that $x_\eps>0, \forall\eps<\eps_0$.
\end{itemize}

 As a corollary we have: if $\tilde{x}$ is a near-standard point
 associated with $x\in\R$ then
\begin{itemize}
\item[(i)] if $x\not= 0$ (resp. $x>0$) then $\tilde{x}$ is invertible
 (resp. strictly positive);
\item[(ii)] if $\tilde{x}\geq 0$ then $x\geq 0$.
\end{itemize}

 Note that the converse implications do not hold as demonstrated by the following
 simple examples: take $\tilde{x}:=[(\eps)_\eps]$. Then $\tilde{x}$ is invertible 
 (moreover strictly positive), but $\tilde{x}\approx 0$. Also,
 $[(-\eps)_\eps]\approx 0$, but $[(-\eps)_\eps]<0$.
 We shall repeatedly make use of the following result:
 
 \blem \label{l:x=0}
 If $\tilx, \tily \in\tilR$ and $|\tilx|\leq \tils\cdot |\tily|$,
 for all invertible $\tils\in\tilR$ with $\ 0<\tils\leq s_0$ for some $0<s_0\in \R$, then
 $\tilx = 0$ in $\tilR$.
 \elem
 \pr
 Suppose $\tilx\not= 0$ in $\tilR$. Then for any representative
 $(x_\eps)_\eps$ of $\tilx$ there is $l\in\R$ and a sequence
 $\eps_k \searrow 0$ such that
 $$
 |x_{\eps_k}|>\eps_k^l, \quad \forall k\in\N.
 $$
 Further, by assumption, for any representatives $(y_\eps)_\eps$
 of $\tilde{y}$ and $(s_\eps)_\eps$ of $\tilde{s}$ and for all
 $a>0$
 $$
 \exists \eps_0\ \forall \eps<\eps_0: \quad
 s_\eps\cdot |y_\eps| - |x_\eps|+\eps^a \geq 0.
 $$
 In particular, this inequality holds for each $\eps_k$ and $k\in\N$ large enough.
 Fix $(y_\eps)_\eps$ and choose $N\in\N$ and $k_0\in\N$ for which
 $|y_{\eps_k}|\leq\eps_k^{-N}, \forall k>k_0$. Then
 $$
 s_{\eps_k}\cdot \eps_k^{-N}>\eps_k^l - \eps_k^a> \eps_k^{l+1},
 $$
 for $a,k$ large enough, hence
 $$
 s_{\eps_k}>\eps_k^{N+l+1}.
 $$
 Since this inequality would have to hold for all invertible $\tils\in\tilR,\ 0<\tils\leq s_0$,
 we arrive at a contradiction (e.g. take $\tils =[(\eps^{N+l+2})_\eps]$).
 \ep

\brem \label{lemmarem} The above proof in fact shows that it suffices to check the assumption for all
$\tilde s = [(\eps^m)_\eps]$ with $m>0$.
\erem

 Turning now to non-constant generalized functions we recall from \cite{OPS-positivity,EbM-AF}
 the definitions of nonnegative and strictly positive generalized
 functions as well as a characterization of positivity:
 $f\in\cG(\Omega)$ ($\Omega\subseteq \R^p$ open) is said to be nonnegative, $f\geq 0$, if $ \forall \,
 (f_\eps)_\eps\ \forall K\subset\subset\Omega\ \forall a>0\ \exists \eps_0\
 \forall \eps<\eps_0: \inf_{x\in K}f_\eps(x) +\eps^a\geq 0$.
 $f$ is said to be strictly positive, $f>0$, if $ \forall
 (f_\eps)_\eps\ \forall K\subset\subset\Omega\ \exists a>0\ \exists
 \eps_0\ \forall \eps<\eps_0: \inf_{x\in K}f_\eps(x) >\eps^{a}$.
 The latter can be characterized as follows: a generalized function
 $f\in\cG(\Omega)$ is strictly positive if and only if $\forall (f_\eps)_\eps\
 \forall K\subset\subset\Omega\ \exists \eps_0\ \forall x\in K\ \forall
 \eps<\eps_0: f_\eps(x)>0$. Moreover, by \cite[Theorem 1.2.5]{book}, $f$ is invertible
 in $\gs(\Omega)$ if and only if $|f|$ is strictly positive.

 We conclude this brief overview by recalling the notion of positive
 (semi)definiteness of symmetric bilinear forms on $\tilR^p$. For this purpose we
 first need the definition of a free generalized vector (cf.\
 \cite{EbM-AF}): let $\tilde{\xi}\in\tilR^p$ and $\tilde\xi\not=0$ in $\tilR^p$.
 We call $\tilde\xi$ free if whenever $\tilde\lambda\in\tilde{\R}$ and $\tilde\lambda\cdot\tilde\xi=0$
 then it follows that $\tilde\lambda=0$ in $\tilR$. Let $A\in\tilR^{p^2}$ be a
 symmetric matrix (i.e. $A=A^{\sf T}$ in $\tilR^{p^2}$). $A$ is positive
 semidefinite if for any $\tilde{\xi}\in\tilR^p$, $\tilde{\xi}^{\sf T} A
 \tilde{\xi}\geq 0$ in $\tilR$. $A$ is positive definite if for all
 free vectors $\tilde{\xi}\in\tilR^p$, $\tilde{\xi}^{\sf T} A \tilde{\xi}$ is strictly 
 positive (i.e., invertible and positive) in $\tilR$. Equivalent characterizations 
 for positive definiteness, resp.\
 semidefiniteness of a symmetric matrix $A\in\tilR^{p^2}$ are:
 
\begin{itemize}
\item[(i)] all eigenvalues of $A$ are strictly positive (resp.
 nonnegative);
\item[(ii)] $A$ has a representative $(A_\eps)_\eps$ such that each $A_\eps$ is symmetric and 
 positive definite (resp. positive semidefinite);
\item[(iii)] for each symmetric representative $(A_\eps)_\eps$ of $A$
 there exists $\eps_0$ such that $A_\eps$ is positive definite (resp.
 positive semidefinite), for all $\eps<\eps_0$.
\end{itemize}

 \section{Differential calculus in algebras of generalized functions} \label{diffcalc}

Basic necessary and sufficient conditions for extremals in the calculus of variations
are usually based on the following elementary results on real valued functions: Let
 $f$ map $\Omega$ to $\R$, where $\Omega$ is an open
 subset of $\R^p$. Let $x_0\in\Omega$ be a local minimum (or maximum)
 of $f\in \cC^1(\Omega)$. Then all first order partial derivatives of
 $f$ vanish at $x_0$. Such points are called critical (or stationary) points.
 When $f\in\cC^2(\Omega)$ then, in addition, the Hessian matrix of
 $f$ is positive (or negative) semidefinite at $x_0$. A sufficient
 condition for a local extremum of $f\in\cC^2(\Omega)$ is as follows:
 if $x_0\in\Omega$ is a critical point of $f$ such that $D^2f(x_0)$ is
 positive (or negative) definite on $\R^p$, then $x_0$ is a local
 minimum (or maximum) of $f$.

 In this section we want to derive analogous criteria for extremals of
 generalized functions. The main difficulty in carrying out this task
 lies in the algebraic structure of $\tilde \R$: the existence of zero 
 divisors necessitates certain adaptations in the classical arguments and
 provides for a more complex structure of the attainable results.
 In this paper we are interested mainly in the extrema of generalized
 functions in classical points (i.e. in $x\in\R^p$) as this suffices
 for the applications in the calculus of variations we aim at. We will however add
 remarks which give analogous criteria for extrema in generalized points. 
 
 As in the classical case, we define a local minimum (resp.\ maximum)
 of $f\in\cG(\Omega)$ ($\Omega$ is an open subset
 of $\R^p$) to be an $x_0\in \Omega$ with the property that
 there exists a neighborhood $\Omega'\subseteq\Omega$ of $x_0$
 such that $f(x_0)\leq f(\tilx)$ (resp.\ $f(\tilde{x})\leq f(x_0)$),
 for all $\tilx\in\tilde{\Omega}'_c$.


 Our first result provides necessary conditions for a minimum:

 \bprop \label{l:gf-necessary cond}
 Let $f\in\cG(\Omega)$. If $x_0\in\Omega$ is a local minimum of $f$ then
 \bi
 \item[(i)] $Df(x_0)=0$ in $\tilR^p$;
 \item[(ii)] $D^2f(x_0)\in\tilR^{p^2}$ is positive semidefinite.
 \ei
 \eprop
 \pr
 Without loss of generality we may suppose that $x_0=0$.\\
 (i) For any $i\in \{1,\dots,p\}$, let $g_i(\cdot):=f(0,\dots,0,\cdot,0,\dots,0)$ be the
 $i$-th partial function. Then $g_i\in\cG(\Omega_i)$ ($\Omega_i = \{t\in \R \mid (0,\dots,t,\dots,0)\in \Omega\}$ 
 open in $\R$) and $0$ is a local minimum of $g_i$ on $\Omega_i$. 
 We have to show that $g'_i(0) = \pd_i f(0) = 0$ in $\tilR$ for all $i$.
 Let $i\in\{1,\dots,p\}$ and assume that $0$ minimizes
 $g_i$ over $(-s_0,s_0)^{\widetilde{}}_c$ ($s_0>0$).
 Let $\tilde{s}=[(s_\eps)_\eps] \in (-s_0,s_0)^{\widetilde{}}_c$
 and take any representative $(g_{i\eps})_\eps$ of $g_i$. Then for all
 $\eps$ we have a Taylor expansion of $g_{i\eps}$ at $0$:
 $$
 g_{i\eps}(s_\eps)=g_{i\eps}(0)+s_\eps g'_{i\eps}(0) + \frac{s_\eps^2}{2}
 g''_{i\eps} (\theta_\eps s_\eps),
 $$
 with $0<\theta_\eps< 1$.
 $g_i(0)\leq g_i(\tilde{s})$ says that
 $$
 \forall a>0\ \exists \eps_0\ \forall \eps<\eps_0:\
 g_{i\eps}(s_\eps) - g_{i\eps}(0) +\eps^a \geq 0.
 $$
 Set $M_\eps:= \sup_{x\in [-s_0,s_0]}|g''_{i\eps} (x)|$.
 Then
 \beq \label{eq:wlog}
 0\leq s_\eps g'_{i\eps}(0) + \frac{s_\eps^2}{2} g''_{i\eps}
 (\theta_\eps s_\eps) +\eps^a
 \leq s_\eps g'_{i\eps}(0) + \frac{s_\eps^2}{2} M_\eps +\eps^a .
 \eeq
 Suppose $\tils$ is strictly positive (hence invertible).
 Then for any representative $(s_\eps)_\eps$ of $\tils$ there exist $m\in\N$ and
 $\eps_1$ such that $s_\eps\geq\eps^m$, for all $\eps<\eps_1$.
 Divide (\ref{eq:wlog}) by $s_\eps$ for $\eps<\min (\eps_0,\eps_1)$:
 $$
 0 
 \leq g'_{i\eps}(0) + \frac{M_\eps}{2} |s_\eps| + \eps^{a'},
 $$
 where $a'=a-m$.
 Otherwise, if $[(s_\eps)_\eps]$ is strictly negative 
 we similarly obtain
 $$
 0 
 \geq g'_{i\eps}(0) - \frac{M_\eps}{2} |s_\eps| -\eps^{a'}\,.
 $$
 Collecting both cases yields
 $$
 |g_i'(0)| \leq \frac{M}{2}|\tilde{s}|,
 $$
 for all invertible generalized numbers $\tilde{s}=[(s_\eps)_\eps]\in
 (-s_0,s_0)^{\widetilde{}}_c$, with $M=[(M_\eps)_\eps]$. The claim now
 follows from Lemma \ref{l:x=0}. \\
 (ii) Taylor expansion of $f$ around $0$ gives:
 $$
 f_\eps(x_\eps)= 
 f_\eps(0)+ Df_\eps(0)(x_\eps)+ \frac12 x_\eps^{\sf T}\cdot
 D^2f_\eps(0)\cdot x_\eps + \frac{1}{3!}((x_\eps
 D)^3f_\eps)(\theta_\eps x_\eps).
 $$
 Now for any $a>0$, $f_\eps(x_\eps)-f_\eps(0)+\eps^a\geq 0$
 for $\eps$ sufficiently small. Also by (i), $Df_\eps(0)(x_\eps)=O(\eps^a)$,
 for all $a>0$. Moderateness of $D^3 f$ yields the existence of some $N>0$ such
 that $\sup_{|y|\le 1} \|D^3 f(y)\| = O(\eps^{-N})$.
 Thus if we suppose that $|x_\eps| = \eps^s$ with $s = \frac{a+N}{3}$ then
 we obtain that $ x_\eps^{\sf T}\cdot D^2f_\eps(0)\cdot x_\eps + \eps^{a-1} \ge 0$
 for any given $a$ if $\eps$ is sufficiently small. 

 $D^2f(0)$ possesses a symmetric representative, so 
 for each $\eps$
 there exists an orthogonal matrix $U_\eps$ such that with $U=[(U_\eps)_\eps]\in \tilde \R^{p^2}$,
 $U\!\cdot \!D^2f(0)\!\cdot \!U^{\sf T}={\mbox{diag}(\lambda_1, \dots,\lambda_p)}$
 with $\lambda_1, \dots,\lambda_p$ the eigenvalues of $D^2f(0)$ in $\tilde \R$ 
 (cf.\ \cite[Lemma 1.6]{EbM-AF}).
 Denote by $\tilx_{\lambda_i}$ the corresponding unit eigenvectors
 (the columns of $U$)
 and set $y_{i\eps}:=\eps^s\cdot x_{\lambda_i\eps}$. Then by the above, for each $a'>0$
 and $\eps$ small
 $$
 0\leq \eps^{a'} + y_{i\eps}^{\sf T}\cdot D^2f_\eps(0)\cdot y_{i\eps}
 = \eps^{a'} + \eps^{2s}\lambda_{i\eps}\,,
 $$
 implying that all
 $\lambda_i$ are positive in $\tilR$. As mentioned in the introduction,
 this is equivalent to positive semidefiniteness of $D^2f(0)$.
 \ep
 
 \brem \label{genpoint}
 ({\em Extrema in generalized points}). In order to extend the validity of \ref{l:gf-necessary cond} to $\tilde x_0\in \tilde \Omega_c$ we first need a concept of local
 minimum in this setting. We employ the notion of support of a generalized point $\tilde x_0 \in \tilde \Omega_c$
 introduced in \cite{garetto_duals}: the support $\supp(\tilde x_0)$ of $\tilde x_0$ is defined as the set of accumulation points
 of any representative $(x_{0\eps})_\eps$. We call $\tilde x_0$ a local minimum of $f\in \gs(\Omega)$ if there exists 
 a neighborhood $\Omega'$ of $\supp(x_0)$ such that $f(\tilde x_0)\le f(\tilde x)$ for all $\tilde x \in \tilde \Omega'_c$. The proof of \ref{l:gf-necessary cond} (i) can be adapted to this situation: set
 $g_{i\eps}(t):= f_\eps(x_{0\eps}^1,\dots,x_{0\eps}^i+t,\dots,x_{0\eps}^n)$. Then although there need not
be a finite interval on which all $g_{i\eps}$ are simultaneously defined, the above arguments may still be
applied for $\eps$ small (using \ref{lemmarem}) and each $i$, showing that $Df(\tilde x_0)=0$. Similarly we may modify the proof
of \ref{l:gf-necessary cond} (ii).

 \erem

 The above result demonstrates that the typical necessary conditions
 for local minima are entirely analogous to those for minima of classical
 functions: the gradient vanishes and the Hessian matrix is positive
 semidefinite. However, as the following example shows  one should not
 expect such a direct analogy with the smooth setting in the case of 
 sufficient conditions: We construct a generalized function whose first derivative
 vanishes at zero and whose second derivative at zero is
 strictly positive but which does not have a minimum in $0$.

 \bex
 Let $\vphi \in\cC^\infty_0(\R)$, $\supp\,\vphi \subseteq [-2,2]$,
 $\vphi|_{(-1/2,1/2)}=x^2$ and $\vphi(-1)=\vphi(1)=-1$. Define
 $f_\eps(\cdot):=\vphi\big(\frac{\cdot}{\eps}\big)$ and
 $f:=[(f_\eps)_\eps]$. Clearly, $f\in \cG(\R),\ f_\eps(0)=0,\
 f'_\eps(0)=0$ and $f''_\eps(0)=\frac{2}{\eps^2}$, so $f''(0)$ is strictly 
 positive (even $f''_\eps(0)\to\infty$). 
 But zero is not a local minimum since the value of $f$ in the generalized 
 point $[(\eps)_\eps]$ is $-1$. 
 \eex
 
 We therefore have to alter the condition of positive
 (semi)definiteness in zero to also take into account neighboring 
 generalized points. Of course we would like to include as few
 points as possible. The following example shows that merely 
 considering near-standard points associated to $0$ is insufficient:

 \bex
 Let $\vphi\in\cC^\infty_0(\R)$, $\supp\,\vphi \subseteq (-1,1)$ and
 $\vphi|_{[-1/2,1/2]}=-1$. Define
 $$
 F_\eps(x):=\sum_{n=1}^{\infty}\eps^n \vphi\Big(\frac{x-1/n}{\eps}\Big),
 \quad \eps\in(0,1],
 $$
 and $F:=[(F_\eps)_\eps]$. Then $F\in\cG(\R), F(0)=0$ and $F'(0)=0$ in $\tilR$.
 Moreover, $F''(\tilde{x})\geq 0$ in $\tilR$, for all near-standard points
 $\tilde{x}\approx 0$ (in fact $F''(\tilde{x})=0$ in $\tilR$ for all $\tilx\approx 0$).
 But $F$ does not attain a local minimum at zero. 
 For any $n_0\in \N$,  $F_\eps(\frac{1}{n_0})=-\eps^{n_0}$ for 
 $2\eps< 1/n_0-1/(n_0+1)$, hence $F(\frac{1}{n_0})<0$ in $\tilR$.  
 Therefore, $0$ is not a local minimum of $f$.
 \eex
 As suggested by the previous examples we therefore give a sufficient
 condition for a local minimum which is based on the behavior of $f$
 in a neighborhood of the critical point.

 \bprop \label{l:gf-sufficient cond}
 Let $f\in\cG(\Omega)$, $x_0\in\Omega$, $Df(x_0)=0$ and let
 $D^2f(\tilde{x})$ be positive
 semidefinite for all $\tilde{x}\in \tilde{\Omega}'_c$ (with $\Omega'$  a
 star-shaped neighborhood of $x_0$ in $\Omega$). Then $x_0$
 is a minimum of $f$ in $\tilde{\Omega}'_c$.
 If in addition $D^2f(\tilx)$ is positive definite for all
 $\tilde{x}\in \tilde{\Omega}'_c$, then the minimum is unique on
 $\tilde{\Omega}'_c$.
 \eprop
 \pr
 Again, taking $x_0=0$ is no restriction.
 Let $\tilde{x}=[(x_\eps)_\eps]\in\tilde{\Omega}'_c$. Take a
 representative $(f_\eps)_\eps$ of $f$. For each $\eps$ we have
 $$
 f_\eps(x_\eps) = f_\eps(0)+ Df_\eps(0)(x_\eps)+ \frac{1}{2}((x_\eps
 D)^2f_\eps)(\theta_\eps x_\eps)
 $$
 with $\theta_\eps\in (0,1)$.
 Since $Df(0)=0$ in $\tilR$ it follows that for all $m>0$ and $\eps$
 small, $|Df_\eps(0)(x_\eps)|\leq \eps^m$. On the other hand, by the assumption that
 $D^2f|_{\tilde{\Omega}'_c}$ is positive semidefinite we
 know that $\forall a'>0\ \exists \eps_0$ such that $1/2((x_\eps
 D)^2f_\eps)(\theta_\eps x_\eps) + \eps^{a'}\geq 0, \forall \eps <\eps_0$.
 Let $a>0$ be given and set $a'$ and $m$ to be $a+1$. Then
 \beas
 0 & \leq & \eps^{a+1} + f_\eps(x_\eps) - f_\eps(0)- Df_\eps(0)(x_\eps)
 \\
  & \leq & \eps^{a+1}+ f_\eps(x_\eps) - f_\eps(0)+ |Df_\eps(0)(x_\eps)| \\
  & \leq & \eps^{a+1}+ f_\eps(x_\eps)-f_\eps(0)+ \eps^{a+1} \\
  & \leq & f_\eps(x_\eps)-f_\eps(0) + \eps^a,
 \eeas
 when $\eps\to 0$. Hence $f(\tilde{x})-f(0)\geq  0$
 and the first part of lemma is proved.
 \\
 To show uniqueness, let $\tilx\in\tilde{\Omega}'_c$ and
 $\tilx\not=0$ in $\tilR^p$. Then there exist $l$ and
 $\eps_k\searrow 0$ such that $|x_{\eps_k}|>\eps^l_k$, for all $k\in\N$.
 Choose a representative
 $(f_\eps)_\eps$ of $f$ with $Df_\eps(0)=0$ for all $\eps$ and write
 $f_\eps(x_\eps)=f_\eps(0)+\frac12 ((x_\eps D)^2f_\eps)(\theta_\eps x_\eps)$.
 In order to bring to bear the positive definiteness of $D^2f$ we embed the
 sequence $(x_{\eps_k})_k$ into a free vector in $\tilde \Omega'_c$ as follows:
  Let $\tily$ be any free vector in $\tilde{\Omega}'_c$ (e.g., $\tily=[(\eps)_\eps]$) and
 define $\tilde x'$ as $[(x_{\eps}')_\eps]$, where
 $$
 x_{\eps}':=\left\{\ba{ll}
 x_{\eps_k}, & \eps=\eps_k\\
 y_\eps, & \mbox{otherwise}
 \ea
 \right.
 $$
 Then $\tilx'$ is free, so the positive definiteness of $D^2f$ in $\tilde z:=[(\theta_\eps x_\eps)_\eps]
 \in \tilde\Omega_c'$ implies strict positivity of $(\tilde x')^{\sf T} \cdot D^2f(\tilde z)\cdot \tilde x'$. 
 Thus there exists some $m>0$ such that for large enough $k\in\N$
 $$
 f_{\eps_k}(x_{\eps_k})-f_{\eps_k}(0)=\frac12 ((x_{\eps_k} D)^2f_{\eps_k})
 (\theta_{\eps_k} x_{\eps_k}) >\eps_k^m\,.
 $$

 This implies that $f(\tilx)\not=f(0)$ in $\tilR$, for all $\tilx\in\tilde{\Omega}'_c$, and uniqueness follows.
 \ep

 \brem As in \ref{genpoint} the above result may be extended to generalized points $\tilde x_0$, if we 
 suppose $\Omega'$ to be a convex open neighborhood of $\supp (\tilde x_0)$.
 \erem
 The uniqueness condition in the previous result raises the question whether positive definiteness of the
 second derivative could be replaced by the assumption of
 $D^2f$ being positive semidefinite but different from zero
 in $\tilR$. These two notions are different in $\tilR$, due to the
 existence of zero divisors. The following example, however, shows that the latter condition is too weak.

 \bex
 For any zero divisor $\alpha\ge 0$ in $\tilde \R$ the assignment
 $$
 f(x):=\alpha x^2,
 $$
 defines an element of $\cG(\R)$ which attains its (even global) minimum $0$
 in the point $0$. Moreover, $f'(0)=0$ in $\tilR$ and $f'' = 2\alpha$ is positive 
 and different from zero everywhere in $\tilR$. Nevertheless the minimum at $0$ is not
 unique: take any $\tilde{s}\not= 0$ in
 $\tilR_c$ with $\alpha\cdot\tilde{s}=0$, then also $f(\tilde{s})=0$.
 \eex

 We conclude this section with a characterization of nonnegativity,
 resp.\ strict positivity, of generalized functions by their
 values in all near-standard points. As a consequence we note that 
 the respective assumptions in the previous results can equivalently formulated
 in terms of this proper subset of generalized points.

 \blem \label{l:positivity near-standard p}
 Let $f\in\cG(\Omega)$. Then $f\geq 0$ in $\cG(\Omega)$ (resp.
 $f>0$ in $\cG(\Omega)$) if and only if $f(\tilx)\geq 0$ in
 $\tilR$ (resp. $f(\tilx)>0$ in $\tilR$) for all near-standard
 points $\tilx\in\tilde{\Omega}_c$.
 \elem
 \pr
 We give a proof for the characterization of nonnegativity of $f$, the
 case of strict positivity can be established along the same lines.
 The condition is clearly necessary. Conversely, 
 suppose that $f$ is not nonnegative in $\cG(\Omega)$. Then
 $$
 \exists (f_\eps)_\eps\ \exists a>0\ \exists
 K\subset\subset \Omega\ \forall k\in\N\ \exists \eps_k<\frac1k\
 \exists x_k\in K:\; f_{\eps_k}(x_k) +\eps_k^a < 0.
 $$
 Without loss of generality we may assume that $\eps_{k+1}<\eps_k$
 and that the sequence $\{x_k\}_{k\in\N}$ converges to some
 $x_0\in K$. Set $x_\eps := x_k$ for $\eps\in (\eps_{k+1},\eps_k]$
 and $\tilx := [(x_\eps)_\eps]$. Then $\tilx$ is a
 near-standard point and $f_{\eps_k} (x_{\eps_k})
 +\eps_k^a <0\ \forall k\in\N$,
 hence $f(\tilx)\not\geq 0$,  a contradiction.
 \ep

 Inspection of the above proof shows that in fact the subset of so-called sequential 
 near-standard points suffices to characterize nonnegativity
 resp.\ positivity of generalized functions. Here, a near-standard point $\tilde{x}\approx x$ is called sequential if it has a
 representative $(x_\eps)_\eps$ of the following form: $\exists  \eps_k\searrow0\
 \exists x_{\eps_k}, x_{\eps_k}\to x$, such that $x_\eps=x_{\eps_k}, \forall
 \eps\in (\eps_{k+1},\eps_{k}]$.

 \section{First and second variation} \label{firstsec}

Let $\Omega$ be an open subset of $\R^n$ and let $D= \Omega$ or $D=\bar\Omega$.
Let $\cA$ be a subset of
$\gs(D)$ and let $\cL: \cA \to \tilde \R$. We want to determine
extremals of $\cL$ under certain admissible variations. To make this 
problem accessible to an analytical treatment we make the following 
assumptions (which will automatically be satisfied in the setting 
most relevant to us, namely (\ref{eq:vp pq-dim}) below). Let $\cA_0$,
the space of {\em admissible variations} be an $\tilde \R$-submodule of 
$\gs(D)$ such that $\cA+\cA_0 \subseteq \cA$. Moreover, we
suppose that for all $u\in \cA$ and all $v\in \cA_0$, the mapping
$$
\cL^{u,v}: s \mapsto \cL(u+sv)
$$
is an element of $\gs(J^{u,v})$, where $J^{u,v}$ is some open 
interval around $0$ in $\R$.

We call $u\in \cA$ a minimizer of the functional $\cL$ in $\cA$
if $\cL(u) \le \cL(w)$ in $\tilde\R$ for all $w\in \cA$. Finally,
we define the first and second variation of $\cL$ at $u\in \cA$ 
in the direction $v\in \cA_0$ as
\beq \label{eq:1st&2nd var der}
\delta\cL(u;v):=\left.\frac{d}{ds}\right|_{s=0}\cL^{u,v}(s);\quad
\delta^2\cL(u;v):=\left.\frac{d^2}{ds^2}\right|_{s=0}\cL^{u,v}(s)
\eeq
Under these assumptions we have
\bprop \label{prop:cv-necessary cond}
Let $u\in\cA$ be a minimizer of the functional $\cL$ in $\cA$. 
Then
\bi
\item[(i)] $\delta\cL(u;v)=0$ in $\tilR$ for all $v\in\cA_0$;
\item[(ii)] $\delta^2\cL(u;v)\geq 0$ in $\tilR$ for all $v\in\cA_0$.
\ei
\eprop
\pr This is immediate from Proposition \ref{l:gf-necessary cond}.
\ep

Turning now to our main object of study, let $\Omega$ be an open, 
connected and bounded subset of $\R^p$
 with smooth boundary $\pd\Omega$. Let $\cL: \cG(\bar{\Omega})^q \to
 \tilR$ be a functional of the form
 \beq \label{eq:vp pq-dim}
 \cL(u) := \int_\Omega L(x, \un(x))\; dx,
 \eeq
 where $L\in\cG_\tau((\bar{\Omega}\times\R^q)^{(n)})$. 
 We will write $\cL_\eps$ for the functional $u\mapsto  \int_\Omega L_\eps(x,u)\; dx$,
 where $u\in \cC^\infty(\bar \Omega)$ and $(L_\eps)_\eps$ is any representative
 of $L$.

 We seek a  minimizer $u$ among all
 functions in $\cG(\bar{\Omega})^q$ which additionally satisfy the 
 boundary condition $u|_{\pd \Omega} =u_0$, $u_0\in\cG(\pd\Omega)^q$
 (here $\partial \Omega$ is viewed
 as a $(p-1)$-dimensional submanifold of $\R^p$).
 The integrand $L$ is called the Lagrangian of the variational
 problem. 
 In the notation introduced above we have
 \begin{equation}\label{*}
 \cA=\{w\in\cG(\bar{\Omega})^q |\, w|_{\pd\Omega}=u_0\}
 \end{equation}
 for some $u_0 \in \gs(\partial \Omega)^q$, with the corresponding space of
 admissible variations given by
 \begin{equation}\label{**}
 \cA_0=\{w\in\cG(\bar{\Omega})^q |\, w|_{\pd\Omega}=0\}\,.
 \end{equation}
 Note that the integral in (\ref{eq:vp pq-dim}) is well-defined since
 the composition of a tempered generalized function with an element of $\gs$ 
 is a well-defined element of $\gs$ (cf.\ \cite[Proposition  1.2.30]{book}). 
 Also,  (\ref{eq:vp pq-dim}) is the most general form of the
 variational problem, involving $p$ independent, $q$ dependent
 variables and all partial derivatives of dependent variables up
 to order $n$. In some cases, in order to
 simplify notation we are going to consider only the case of
 one unknown function; the generalization to two or more dimensions is then
 a technical matter.

 \brem \label{rem:c-bounded}
 In the formulation of the variational problem, instead of tempered
 generalized functions one may alternatively use $c$-bounded generalized functions:
 assume that $L\in\cG((\bar{\Omega}\times\R^q)^{(n)})$. 
 Since c-bounded generalized functions may be composed unrestrictedly (see \cite[Proposition 1.2.8]{book}), 
 a minimizer then can be sought among all $c$-bounded generalized functions with
 prescribed boundary conditions. In fact, for the integrand to be well-defined not only $u$ but also
 all partial derivatives of $u$ up to order $n$ have to be $c$-bounded, i.e.\ the
 admissible set is $\cA:=\{w|\, w^{(n)}\in \cG[\bar{\Omega},(\R^q)^{(n)}]
 \wedge w|_{\pd\Omega}=u_0\}$. Although convenient for an invariant geometric formulation of
 the theory, this latter requirement necessitates certain unwanted restrictions on the
 degree of singularity of $u$ (e.g., $\delta$-like singularities are excluded).
 We will therefore mainly operate in the tempered setting, but will simultaneously 
 keep track of the $c$-bounded case.
 \erem

 In order to derive Euler-Lagrange equations in our setting, as in the classical case we first 
 derive a fundamental lemma of the calculus of variations:

 \blem \label{l:fundamental pq-dim}
 (Fundamental lemma of the calculus of variations)
 Let $u\in\cG(\Omega)$. If
 $$
 \int_\Omega u(x)\varphi(x)\, dx=0 \quad \mbox{in } \tilR
 $$
 for all $\varphi\in\cG_c(\Omega)$ 
 then $u \equiv 0$ in $\cG(\Omega)$. 
 \elem
 \pr
 Suppose to the contrary that $u=[(u_\eps)_\eps]$ with $(u_\eps)_\eps \not\in\cN(\Omega)$. 
 Then $\exists K\subset\subset \Omega\
 \exists l>0\ \forall k\in\N\ \exists \eps_k \searrow 0\
 \exists x_k\in K$ such that $|u_{\eps_k}(x_{k})| > \eps_k^l$.
 On the other hand, $(u_\eps)_\eps$ is moderate, so $\exists N>0$ such that
 $\sup_{x\in K}|Du_\eps(x)| \leq \eps^{-N}$. 
 We may choose some $k_0$ such that for $k\ge k_0$ and any $x\in \Omega$ with 
 $|x-x_k|\leq \frac12 \eps_k^{l+N}$, the entire line segment connecting
 $x$ and $x_k$ is contained in $\Omega$. Then
 $$
 u_{\eps_k}(x) = u_{\eps_k}(x_k) +
 Du_{\eps_k}(x_k+\theta_k(x-x_k))(x-x_k),
 $$
 for some $\theta_k\in (0,1)$, so that
 $|u_{\eps_k}(x)|\geq \eps_k^l - \eps_k^{-N}|x-x_k|
 \geq \frac12 \eps_k^l$ for all $x$ with $|x-x_k|\leq \frac12 \eps_k^{l+N}$.

 Choose $\varphi\in\cC^\infty_0(B_{\frac12}(0)), \varphi\geq 0$ and
 $\varphi\equiv 1$ on $B_{\frac14}(0)$ ($B_{r}(0)$ denoting the open
 ball of radius $r$ around $0$). Define $\varphi_\eps\in \cD(\Omega)$ by 
 $$
 \varphi_\eps(x):= \eps^{-l-N}
 \varphi(\eps^{-l-N}(x-x_\eps)) \cdot \sign (u_\eps(x_\eps))
 $$
 where $x_\eps = x_{\eps_k}$ for all $\eps\in (\eps_{k+1},\eps_k]$ and all 
 $k\ge k_0$. Then $\varphi:=[(\varphi_\eps)_\eps]\in \cG_c(\Omega)$.
 Since $u_{\eps_k}$ does not change sign on the support of $\varphi_{\eps_k}$
 and since
 $u_{\eps_k}(x)\cdot \varphi_{\eps_k}(x) \geq \frac12
 \eps_k^{-N}$ whenever $|x-x_k|\leq \frac14 \eps_k^{l+N}$, we obtain
 \beas
 \Big|\int_\Omega u_{\eps_k}(y) \varphi_{\eps_k}(y)\; dy\Big| 
 &=& \int_{B_{\frac12\eps_k^{l+N}}(x_k)}
 u_{\eps_k}(y) \varphi_{\eps_k}(y)\; dy\\
 &\ge& \int_{B_{\frac14\eps_k^{l+N}}(x_k)}
 u_{\eps_k}(y) \varphi_{\eps_k}(y)\; dy\\
 &\geq& \mathrm{const}\cdot\eps_k^{p(l+N)-N},
 \eeas
 which implies that $(\int_\Omega u_\eps(y)\varphi_\eps(y)\,dy)_\eps \not\in \cN(\Omega)$,
 a contradiction. Therefore, $u\equiv 0$ in $\cG(\Omega)$.
 \ep

 As in the smooth setting, this result allows to derive the Euler Lagrange equations
 as necessary conditions for a minimizer of (\ref{eq:vp pq-dim}):

 \bthm \label{eulerlagrange}
 Let $u\in\cA$ be a minimizer of (\ref{eq:vp
 pq-dim}). Then $u$ is a
 solution of the system of Euler-Lagrange equations
 \beq \label{eq:EL pq-dim}
 E_\alpha(L)=0,\ \alpha=1,\dots,q, \quad \mbox{ in }
 \cG(\Omega),
 \eeq
 where $E_\alpha$ is the $\alpha$-th Euler operator
 $$
 E_\alpha=\sum_{J} (-D)_J \frac{\pd}{\pd u^\alpha_J}
 $$
 (the sum extends over all multi-indices $J=(j_1,\dots,j_k),
 1\leq j_k \leq p, 1\leq k\leq n$).
 \ethm
 \pr
 By Proposition \ref{prop:cv-necessary cond} (i), 
 $\delta \cL(u;v)=0, \forall v\in \cA_0$. Now
 \beas
 \delta \cL(u;v) 
 &=& \int_\Omega \frac{d}{ds}{\Big\vert}_{_{s=0}} L(x,
 (u+sv)^{(n)}(x))\; dx \\
 &=& \int_\Omega \sum_{\alpha =1}^{q} \sum_{J} \frac{\pd L}{\pd u^\alpha_J}
 (x,\un(x)) \cdot \pd_J v^\alpha(x) \; dx\,.
 \eeas
 Integration by parts gives (with $D$ denoting the total derivative)
 \beas \delta \cL(u;v)
 & = & \int_\Omega \sum_{\alpha =1}^{q} \Big[ \sum_{J} (-D)_J
 \frac{\pd L}{\pd u^\alpha_{J}} (x,\un(x))\Big]
 \cdot v^\alpha(x) \; dx\,
 \eeas
 The claim therefore follows from Lemma \ref{l:fundamental pq-dim}.
 \ep

 \brem
 In the case $n=1$, i.e.\ when $\cL(u)=\int_\Omega L(x,\prfirst u(x))\,
 dx$, the $\alpha$-th Euler operator takes the well known form:
 $$
 \ds E_\alpha(L)=\frac{\pd L}{\pd u^{\alpha}} -
 \sum_{i=1}^p \frac{\pd}{\pd x_i} \frac{\pd L}{\pd
 u^{\alpha}_i}
 $$
 where $u^{\alpha}_i = \pd u^{\alpha}/\pd x_i$.
 \erem

 The following result provides sufficient conditions for minimizers of the variational problem
 (\ref{eq:vp pq-dim}).
 \bprop \label{prop:cv-sufficient cond}
 With $\cL$ as in (\ref{eq:vp pq-dim}) and $\cA$, $\cA_0$ as in (\ref{*}), (\ref{**}), 
 suppose that 
 \begin{itemize}
 \item[(i)] $\delta\cL(u;v)=0$ in $\tilR$ for all $v\in\cA_0$.
 \item[(ii)]  $\delta^2\cL(w;v)\geq 0$ for all $w\in\cA$ and $v\in\cA_0$. 
  \end{itemize} 
 Then $u$ is a minimizer for $\cL$. 
 If, in addition, $s\mapsto \delta^2\cL(u+sv;v)\not= 0$ in $\cG([0,1])$,
 $\forall v\in\cA_0$, $v\not= 0$, then the minimizer is unique.
 \eprop
 \pr
 Given any $w\in\cA$, set $\Psi:= s\mapsto \cL(u+s(w-u)) \in \cG(\R)$. 
 By assumption we have $\Psi'(0)=\delta \cL(u;w-u)=0$ (since $w-u\in\cA_0$)
 and $\Psi''(\tils)=\delta^2\cL(u+\tils(w-u);w-u)\geq 0$, $\forall
 \tils\in\tilR_c$. We now apply Lemma \ref{l:gf-sufficient cond}
 to conclude that $0$ is a minimum of $\Psi$,
 i.e. $\Psi(0)\leq \Psi(\tils), \forall \tils\in \tilR_c$.
 In particular, for $\tils=1$ we have that $\cL(u)\leq \cL(w)$,
 and since $w\in\cA$ was arbitrary the claim follows.\\
 If, in addition, $s\mapsto \delta^2\cL(u+s(w-u);w-u)\not= 0$ in $\cG([0,1])$
 for all $w\in\cA$, i.e. $s\mapsto \Psi''(s)\not= 0$ in $\cG([0,1])$,
 then according to \cite[Lemma 3.1]{HMPV} (which in fact is closely related to our 
 Lemma \ref{l:fundamental pq-dim}) we obtain that
 $$
 \Psi(1)-\Psi(0)=\int_0^1 (1-s)\Psi''(s)\; ds \not= 0\,.
 $$
 Therefore $\cL(w)-\cL(u)\not= 0$ for all $w\not= u$, and uniqueness
 follows.
 \ep

\brem The uniqueness criterion in Proposition \ref{prop:cv-sufficient cond} can be adapted to provide an
alternative to the uniqueness condition in Proposition \ref{l:gf-sufficient cond}, to wit: under the
assumptions given there, suppose additionally that $s\mapsto D^2f(x_0+s(\tilde y -x_0))[\tilde y-x_0,\tilde y-x_0]\not=0$
in $\gs[0,1]$ for all $\tilde y\in \tilde \Omega'_c$. Then the minimum is unique in $\tilde \Omega'_c$.
 
\erem

\bex \label{ex:quadraticforms} (Quadratic forms). Suppose that for each $\eps\in (0,1]$, $a_\eps: \cC^\infty(\bar\Omega)\times 
\cC^\infty(\bar\Omega) \to \R$ is symmetric, bilinear and continuous. 
We additionally suppose that 
$(a_\eps(u_\eps,v_\eps))_\eps \in \cE_M$ for all $(u_\eps)_\eps, (v_\eps)_\eps 
\in \cE_M(\bar\Omega)$, and $(a_\eps(u_\eps,v_\eps))_\eps \in \cN$ for $(u_\eps)_\eps 
\in \cE_M(\bar\Omega)$ and $(v_\eps)_\eps \in \cN(\bar \Omega)$ 
(where $\tilde \R = \cE_M/\cN$), so that
\begin{eqnarray*}
a: \cG(\bar \Omega)\times \cG(\bar\Omega) &\to& \tilde\R \\
a(u,v) &:=& [(a_\eps(u_\eps,v_\eps))_\eps]	
\end{eqnarray*}
is a well-defined bilinear map. 
Similarly, for each $\eps \in (0,1]$ let $f_\eps: \cC^\infty(\bar\Omega) \to \R$ be linear and
continuous and suppose that $(f_\eps(u_\eps))_\eps \in \cE_M$ (resp.\ $\cN$) for
each $(u_\eps)_\eps \in \cE_M(\bar \Omega)$ (resp.\ $\cN(\bar \Omega)$). Then
set $f: \cG(\bar\Omega)\to \tilde \R$, $f(u):=[(f_\eps\circ u_\eps)_\eps]$.

We consider the functional
\begin{equation} \label{quadform}
\begin{array}{rcl}
\cL: \cG(\bar \Omega) &\to& \tilde \R \\[3pt]	
\cL(u) &:=& \frac12 a(u,u) + f(u)\,.	
\end{array}
\end{equation}

Again we will write $\cL_\eps$ for the functional $u\mapsto \frac12 a_\eps(u,u) + f_\eps(u)$.
In this situation, condition (i) from Proposition \ref{prop:cv-necessary cond} reads 
\begin{equation}\label{VE}
\delta \cL(u,v) =  a(u,v) + f(v) = 0
\end{equation}
and $\delta^2\cL(u,v) = a(v,v)$. Thus if $a$ is positive semidefinite in the sense 
that $a(u,u)\ge 0$ for all $u\in \cG(\bar \Omega)$
then any solution of (\ref{VE}) is a minimizer of $\cL$.

As a concrete example let
\begin{equation}\label{eq:quadraticform}
   \cL(u) = \frac{1}{2}\int_\Omega \big(\alpha(x)|\nabla u(x)|^2 + \beta(x)u(x)^2\big)\,dx
                     + \int_\Omega \gamma(x)u(x)\,dx
\end{equation}
with $\cA$, $\cA_0$ as in (\ref{*}), (\ref{**}), 
$\alpha, \beta, \gamma\in \cG(\bar \Omega)$ and $\alpha\geq 0, \beta\geq 0$ on $\bar{\Omega}$. 
Under these assumptions, $a$ is positive semidefinite and by the above,
any solution of (\ref{VE}) is a minimizer of $\cL$. Moreover, the minimizer 
is unique by Prop.~\ref{prop:cv-sufficient cond} in case $\alpha$ is invertible on $\bar{\Omega}$.
In fact, if
$$
\delta^2\cL(u;v) = a(v,v)  
    = \int_\Omega \big(\alpha(x)|\nabla v(x)|^2 + \beta(x)|v(x)|^2\big)\,dx = 0
$$
it follows that $\nabla v \equiv 0$ in $\cG(\Omega)$, so $v$ is a generalized
constant. If $v$ is an element of $\cA_0$, this constant must be zero.
Therefore, $\delta^2\cL(u,v)\neq 0$ for all $v\in \cA_0$, $v\neq 0$; this in turn
implies the hypotheses of Prop. \ref{prop:cv-sufficient cond}.

On the other hand, uniqueness may be lost if $\alpha$ is not invertible: let $\Omega=(-1,1)$, and
$$
\cL(u) = \frac12 \int_{-1}^1 \alpha |u'(x)|^2 \,dx + \int_{-1}^1 \alpha u(x)\,dx
$$
with $\alpha\ge 0$ a zero divisor, say $\alpha \omega = 0$. 
Then $\delta^2 \cL(u,v)=\int_{-1}^1 \alpha |v'(x)|^2\,dx\ge 0$
for all $v$. However, both $u(x)= x^2/2-1/2$ 
and $\bar u(x)= \omega x^2/2 - \omega/2$ are minimizers, because they
satisfy the Euler-Lagrange equation
\[
   \alpha\big(-u''(x) + 1\big) = 0,\quad u(-1) = u(1) = 0.
\]
\eex

We now turn to an analysis of variational problems in the distributional setting. 
In the language of Colombeau algebras,
this amounts to expressing extremal properties in terms of the concept of association:
\bd An element $u$ of $\cG(\Omega)$ is called a minimizer in the sense of association for the 
functional (\ref{eq:vp pq-dim}) if for any (hence all) representatives $(\cL_\eps)_\eps$ of $\cL$ and 
$(u_\eps)_\eps$ of $u$ we have
\begin{equation} \label{assocmini}
\lim_{\eps\to 0} (\cL_\eps(u_\eps + \vphi) - \cL_\eps(u_\eps)) \ge 0   
\qquad \forall \vphi \in \cD(\Omega)\,.
\end{equation}
\ed
For a large class of quadratic forms, the following result provides a necessary and sufficient
condition for a minimizer in the sense of association:
\blem \label{qflemma}
Let $\cL$ be a functional of the form (\ref{quadform}) where $a$ is positive semidefinite 
and such that $A(\vphi,\vphi):=\lim_{\eps\to 0}a_\eps(\vphi,\vphi)$ exists for all $\vphi\in \cD(\Omega)$.
Let $u\in \cG(\Omega)$. Then the following are equivalent:
\begin{itemize}
\item[(i)] $u$ is a minimizer in the sense of association for $\cL$.
\item[(ii)] $\lim_{\eps\to 0} (a_\eps(u_\eps,\vphi)+f_\eps(\vphi)) = 0$ for all $\vphi\in \cD(\Omega)$.
\end{itemize}
\elem
\pr (ii)$\Rightarrow$(i): Let $\vphi\in \cD(\Omega)$. Then
$\cL_\eps(u_\eps + \vphi) - \cL_\eps(u_\eps) = a_\eps(u_\eps,\vphi) +
\frac12 a_\eps(\vphi,\vphi) + f_\eps(\vphi) \to \frac12 A(\vphi,\vphi)\ge 0$ ($\eps\to 0$).

(i)$\Rightarrow$(ii): Again fix $\vphi\in \cD(\Omega)$. Inserting $\tau \vphi$ instead of
$\vphi$ in (\ref{assocmini}), we obtain
$$
\lim_{\eps\to 0}(a_\eps(u_\eps,\vphi)+f_\eps(\vphi)) + \frac{\tau}{2} A(\vphi,\vphi) \ge 0 
\quad \forall \tau>0\,,
$$
and 
$$
\lim_{\eps\to 0}(a_\eps(u_\eps,\vphi)+f_\eps(\vphi)) + \frac{\tau}{2} A(\vphi,\vphi) \le 0 
\quad \forall \tau<0\,,
$$
which together give the claim.
\ep

The following is a classical example due to Weierstrass, cf.\ also \cite{graves}:
\bex 
Let $\Omega=(-1,1)$, $c\not=d\in \R$ and consider the functional
$$
\cL(u) = \int_{-1}^1 x^2 u'(x)^2\,dx\,.
$$
We want to minimize $\cL$ subject to the boundary conditions $u(-1)=c$, $u(1)=d$. Here, $a_\eps$ is
in fact independent of $\eps$. Denote by $v$ the step function equal to $c_1$ on $(-1,0)$ and equal to
$c_2$ on $(0,1)$ ($v$ is the ``expected'' minimizer). Let $\rho$ be a standard mollifier ($\rho\in \cD(\R)$, $\int \rho =1$) and set 
$u_\eps := v*\rho_\eps$. Then 
$$
\lim_{\eps\to 0} (\cL(u_\eps + \vphi) - \cL(u_\eps)) = \lim_{\eps\to 0} 
\int_{-1}^1 (2x^2 u_\eps'(x)\vphi'(x) + x^2\vphi'(x)^2)\,dx \ge 0 \quad \forall \vphi \in \cD(\Omega)\,,
$$
so $u=[(u_\eps)_\eps]$ is a minimizer in the sense of association. Moreover, $\cL(u_\eps) \to 0$ for $\eps\to 0$ and
$u$ is associated to $v$. By Lemma \ref{qflemma}  (or by direct calculation), $u$
solves the corresponding Euler-Lagrange equation in the  sense of association.

Note, however, that $u$ is not a minimizer in $\cG(-1,1)$. Otherwise by Theorem \ref{eulerlagrange} it would have
to satisfy the corresponding Euler-Lagrange equations (with equality in $\cG$). This, however, 
only has the classical
solution, so also in $\cG(-1,1)$, no solution satisfying the boundary conditions exists.
\eex

For more general variational problems of the form (\ref{eq:vp pq-dim}) we need additional 
conditions in order to have Euler-Lagrange equations (as in Lemma \ref{qflemma} (ii)) 
as necessary conditions for minimizers (in the sense of association):

\bprop For the variational problem (\ref{eq:vp pq-dim}), suppose that the total
second order derivative of $L$ w.r.t.\ the variables $u^{(n)}$
is globally bounded (independently of $\eps$).
If $u\in \cG(\Omega)$ 
is a minimizer of $\cL$ in the sense of association, then $u$ satisfies the Euler-Lagrange
equations for $\cL$ in the sense of association:  
$$ 
E_\alpha(L)\approx 0,\ \alpha=1,\dots,q\,,
$$ 
 i.e.,  
$$ 
 \lim_{\eps\to 0} \int_\Omega E_\alpha(L_\eps)(u_\eps)(x) 
\vphi(x)\,dx = 0 \qquad \forall \vphi \in \cD(\Omega)\, 
 \forall \alpha\,.
$$
\eprop 
\pr
To simplify notations we give the proof for $p=q=n=1$, the general case being
entirely analogous. By Taylor expansion we have
$$
L_\eps(x,y+h,z+k) = L_\eps(x,y,z) + \partial_2L_\eps(x,y,z)h +\partial_3L_\eps(x,y,z)k
+R_\eps(x,y,z,h,k) 
$$
where 
$$
R_\eps(x,y,z,h,k) = (\int_0^1 D^2_{y,z}L_\eps(x,y+th,z+tk)(1-t)\,dt)[(h,k),(h,k)]
$$
and by assumption, $|R_\eps(x,y,z,h,k)| \le C\|(h,k)\|^2$ for all $\eps$.
Let $\vphi\in \cD(\Omega)$, $\tau\in \R$. Then
\begin{eqnarray*}
0 &\le& \lim_{\eps\to 0}(\cL_\eps(u_\eps + \tau\vphi) - \cL_\eps(u_\eps))\\
&=& \lim_{\eps\to 0} \left(\tau\int_\Omega \partial_2L_\eps(x,u_\eps(x),u_\eps'(x))\vphi(x)
+  \partial_3L_\eps(x,u_\eps(x),u_\eps'(x))\vphi'(x)\,dx \right.\\
&& \hphantom{\lim_{\eps\to 0} \ (}
\left. +\, \tau^2\, (\int_\Omega R_\eps(x,u_\eps(x),u_\eps'(x),\tau\vphi(x),\tau\vphi'(x)))
[(\vphi(x),\vphi'(x)),(\vphi(x),\vphi'(x))] \,dx\right) \\
\end{eqnarray*}
Since the second integral is uniformly bounded, dividing by $\tau$ and distinguishing 
the cases $\tau>0$ and $\tau<0$ (cf.\ the proof of Lemma \ref{qflemma}) implies the result.
\ep

 \section{Variational symmetries and N\"{o}ther's Theorem} \label{varsymmnoether}

In this section we turn our attention to symmetries of variational
problems involving generalized functions. After a quick recapitulation
of the relevant notations and definitions in the smooth setting
(following \cite[Ch.\ 4]{Olv}), our goal is to introduce an appropriate notion of  
generalized variational symmetry groups, retaining the fundamental 
properties of classical variational symmetries. Based on this, we will
then derive a corresponding version of N\"other's theorem.

 Turning first to the smooth setting, let $\cL: \cC^\infty(\bar{\Omega}) 
 \to \R$ be a functional of the form
 \beq \label{eq:cvp}
 \cL(u) =\int_\Omega L(x, u^{(n)}(x))\, dx,
 \eeq
 where $\Omega$  is an open, connected set in $\R^p$ with smooth boundary 
 $\pd\Omega$ and $L$ is a smooth function on the $n$-jet space
 $(\bar{\Omega}\times\R^q)^{(n)}$.
 A variational symmetry group of (\ref{eq:cvp}) is a local group of
 transformations $G$ acting on $M\subset \Omega\times\R^q$ with the 
 following property: let $\Omega_1$ be any subdomain with
 $\bar{\Omega}_1\subset\Omega$, $f$ a smooth function 
 on $\Omega_1$ whose graph lies in $M$, and suppose that for $g\in G$,
 $f' := g\cdot f$ is a well-defined function on some $\Omega'_1\subset \Omega$, then
 \beq \label{eq:cvsg}
 \int_{\Omega'_1} L(x', \prn f'(x')) \, dx' =
 \int_{\Omega_1} L(x,\prn f(x))\, dx.
 \eeq
 Heuristically, therefore, a variational symmetry group is a 
 transformation group that does not change the variational integral.
 Since in the generalized setting we will exclusively deal with projectable group
 actions, let us investigate (\ref{eq:cvsg}) in the case where $G$ is of the
 form $(x,u) \mapsto (\Xi_\eta(x),\Psi_\eta(x,u))$, $\Xi, \Psi$ smooth,
 $\Xi_\eta$ a (w.l.o.g.\ orientation preserving) diffeomorphism for each $\eta$. 
 Then $x' = \Xi_\eta(x)$ is
 a coordinate transformation and we can rewrite (\ref{eq:cvsg}) as
 \beq \label{eq:cvsg1}
 \int_{\Omega_1} L(\Xi_\eta(x), \prn f'(\Xi_\eta(x))) \det \mathrm{D}\Xi_\eta(x) \, dx =
 \int_{\Omega_1} L(x,\prn f(x))\, dx
 \eeq  
 (with $f' = x' \mapsto \Psi_\eta(\Xi_{-\eta}(x'),f'(\Xi_{-\eta}(x')))$),
 for each $\Omega_1$ as above and $\eta$ sufficiently small.

Our strategy basically will be to turn (\ref{eq:cvsg1}) into a definition
in the generalized setting. Thus we want to allow for $G$ a generalized
projectable group action  $\Phi\in\tilde{\cG}_\tau(\R\times\R^{p+q})^{p+q}$ 
of the form $(\eta,(x,u))$ $\mapsto$ $(\Xi_\eta(x),\Psi_\eta(x,u))$. 
Since contrary to the smooth setting, we cannot
argue by continuity in $\eta$ (in fact, we explicitly want to allow
e.g.\ jumps), in order to make (\ref{eq:cvsg1}) well-defined for 
general $L\in\cG_\tau((\bar{\Omega}\times\R^q)^{(n)})$, we henceforth
assume that $\Xi\in\cG[(-\eta_0,\eta_0)\times\Omega,\Omega]$ for some
$\eta_0>0$. This guarantees that for each $\Omega_1$ with $\bar\Omega_1\subset 
\Omega$ there exists some $\eta_1\in (0,\eta_0)$ such that $\Xi_\eps(\eta,x)$ 
remains in a fixed compact subset of $\Omega$ for all $\eta<\eta_1$ and $\eps$ 
sufficiently small. Therefore, the following definition makes sense:

\bd \label{genvarsymm}
Let $\Phi=(\eta,(x,u))\mapsto (\Xi_\eta(x),\Psi_\eta(x,u))\in
\tilde{\cG}_\tau(\R\times\R^{p+q})^{p+q}$ be a
generalized (projectable) group action with  $\Xi|_{(-\eta_0,\eta_0)\times\Omega} 
\in\cG[(-\eta_0,\eta_0) \times\Omega,\Omega]$ for some $\eta_0>0$. $\Phi$ is called a
variational symmetry group of the functional (\ref{eq:vp pq-dim})
with $\cL \in \cG_\tau((\bar \Omega\times\R^q)^{(n)})$ 
if for every $\Omega_1 \subseteq \bar\Omega_1\subseteq \Omega$ 
and every $f\in \cG(\Omega_1)^q$ we have
$($with $f' = x' \mapsto \Psi_\eta(\Xi_{-\eta}(x'),f'(\Xi_{-\eta}(x')))):$
\begin{equation}\label{genvarsymmeq}
 \int_{\Omega_1} L(\Xi_\eta(x), \prn f'(\Xi_\eta(x))) \det \mathrm{D}\Xi_\eta(x) \, dx =
 \int_{\Omega_1} L(x,\prn f(x))\, dx \quad  \mbox{in } \cG((-\eta_1,\eta_1))
\end{equation}
whenever $\eta_1>0$ is such that $\Xi_{\eps\eta}(\bar\Omega_1)$ remains in
a fixed compact subset of $\Omega$ for $\eta<\eta_1$ and $\eps$ small.

\ed
\brem ({\em Inheritance properties of variational symmetries for smooth functionals}) \\
Consider a smooth variational problem ((\ref{eq:vp pq-dim}) with $L$ smooth)
and a smooth projectable variational symmetry $\Phi$. 
Then since composition of smooth functions and c-bounded generalized
functions is well-defined, (\ref{genvarsymmeq}) holds for 
any $f\in \cG[\Omega_1,\R^q]$,  since it holds on the level of representatives.
The analogous statement is true for $f\in\cG(\bar{\Omega}_1)^q$ if we assume 
that both $L$ and $G$ are slowly increasing, uniformly for
$x$ in compact sets (to make the respective compositions well-defined).
\erem

 \brem
 As was mentioned in Remark \ref{rem:c-bounded}, an alternative admissible
 set over which the variational problem $\cL$ can be considered
 is the set of all $c$-bounded generalized functions from $\bar{\Omega}$
 to $\R^q$ such that all partial derivatives up to order $n$ are
 also $c$-bounded. In this case, in order for the transformed
 function $f'$ to be well-defined 
 we have to suppose that $\Phi^{(n)}\in \cG[\R\times\R^{p+q},
 (\R^{p+q})^{(n)}]$ (as can readily be verified by direct calculation using the chain rule). 
 \erem

Our next aim is to derive an infinitesimal criterion for variational symmetries. In its
proof we will make use of the following lemma.

 \blem \label{l:det J}
 Let $\Xi\in \tilcG_\tau(\R\times\R^p)^p$ be a $\cG$-complete
 group of transformations acting on $\R\times\R^p$ with
 infinitesimal generator $\xi$. Then
 $$
 \frac{d}{d\eta}(\det J(\eta,x))=\det J(\eta,x)\cdot \Div\xi
 (\Xi_\eta(x)),
 $$
 where $J(\eta,x)$ is the Jacobian matrix of $\Xi$ with entries
 $$
 J_{ij}(\eta,x)=\frac{\pd}{\pd x^j}(\Xi^i_\eta(x)) \in
 \tilde \cG_\tau(\R\times\R^p)^p, \, 1\leq i,j\leq p.
 $$
 and $\Div \xi = \sum_{i=1}^p D_i \xi_i$ is the total divergence of $\xi$.
 \elem
 \pr

 Since the proof of the corresponding result in the smooth case given in \cite[Theorem 4.12]{Olv}
 can not directly be recovered in the generalized setting we give a direct computational
 argument. We have
 \beq \label{eq:der of det}
 \frac{d}{d\eta} (\det J(\eta,x)) = \sum_{k=1}^{n}\sum_{\sigma\in S_n}
 \mathrm{sgn}(\sigma) J_{1\sigma(1)}(\eta,x)\cdot\ldots\cdot \frac{d}{d\eta}
 J_{k\sigma(k)}(\eta,x)\cdot\ldots\cdot J_{n\sigma(n)}(\eta,x).
 \eeq
Here,
 $$
 \frac{d}{d\eta} J_{k\sigma(k)}(\eta,x) 
 = \frac{\pd}{\pd x^{\sigma(k)}} (\xi^k(\Xi_\eta(x)))\\
 = \sum_{l=1}^{n} \frac{\pd\xi^k}{\pd x^l}(\Xi_\eta(x))
 J_{l\sigma(k)}(\eta,x).
 $$
 Inserting into (\ref{eq:der of det}) we obtain:
 \beas
 \frac{d}{d\eta} (\det J(\eta,x)) &=& \sum_{k=1}^{n}\sum_{\sigma\in S_n}
 \mathrm{sgn}(\sigma) \bigg(
 J_{1\sigma(1)}(\eta,x)\cdot\ldots\cdot
 J_{k\sigma(k)}(\eta,x)\cdot\ldots\cdot J_{n\sigma(n)}(\eta,x)
 \frac{\pd\xi^k}{\pd x^k}(\Xi_\eta(x))\\
 && \quad + \sum_{l\not= k}
 J_{1\sigma(1)}(\eta,x)\cdot\ldots\cdot
 J_{l\sigma(k)}(\eta,x)\cdot\ldots\cdot J_{n\sigma(n)}(\eta,x)
 \frac{\pd\xi^k}{\pd x^l}(\Xi_\eta(x)) \bigg)\\
 &=& \det J(\eta,x) \cdot \Div\xi(\Xi_\eta(x)) + 0\,.
 \eeas
 \ep

 \bthm (Infinitesimal criterion) \label{th:inf cri}
 Let $\Phi\in \tilcG_\tau(\R\times\R^{p+q})^{p+q}$ be a generalized,
 projectable, strictly $\cG$-complete group action of the form 
 given in Definition \ref{genvarsymm}. Then $\Phi$ is
 a variational symmetry group of the functional (\ref{eq:vp pq-dim}) if and
 only if
 \beq \label{eq:ic-1}
 \prn \vf(L) +L\cdot\Div\xi=0 \quad \mbox{in
 }\cG((\Omega\times\R^q)^{(n)}),
 \eeq
 for the infinitesimal generator
 $$
 \vf=\sum_{i=1}^{p}\xi^i(x)\frac{\pd}{\pd x^i} +
 \sum_{\alpha=1}^{q}\psi_\alpha(x,u) \frac{\pd}{\pd u^\alpha}
 $$
 of $\Phi$ ($\xi^i\in\cG_\tau(\R^p), \psi_\alpha\in\cG_\tau(\R^{p+q})$).
 \ethm
\pr Suppose first that $\Phi$ is a variational symmetry of (\ref{eq:vp pq-dim}), 
and let $f\in \cG(\Omega_1)^q$, where $\Omega_1\subset \bar \Omega_1 \subset \Omega$. 
Then there exists some $\eta_1>0$ 
such that for all $\tilde \eta \in (-\eta_1,\eta_1)^\sim_c$ and all
$\Omega_2 \subseteq \Omega_1$, by (\ref{genvarsymmeq}) we have
\beq \label{eq:ic-2}
\int_{\Omega_2} L(\Xi_{\tilde \eta}(x), \prn f'(\Xi_{\tilde \eta}(x))) 
\det \mathrm{D}\Xi_{\tilde\eta}(x) \, dx =
\int_{\Omega_2} L(x,\prn f(x))\, dx \quad  \mbox{in } \tilde \R\,.
\eeq

 Since this holds for all open subsets $\Omega_2$ of $\Omega_1$, 
 we may apply \cite[Theorem 1]{PSV} to conclude that for each $\tilde \eta$
 as above the integrands in (\ref{eq:ic-2}) agree in
 $\cG((\Omega_1\times\R^q)^{(n)})$, i.e.,
 \beq \label{eq:ic-3}
 L(\Xi_{\tilde \eta}(x), \prn f'(\Xi_{\tilde \eta}(x)))\det (J_{\tilde\eta}(x)) =
 L(x, \prn f(x)) \quad \mbox{in }
 \cG((\Omega_1\times\R^q)^{(n)}).
 \eeq
Suppose now that $(\tilde x, \tilde u^{(n)})$ is any compactly supported generalized
point in $(\Omega_1\times \R^q)^{(n)}$. 
Then we may choose for the components of $f\in \cG(\Omega_1)^q$ 
suitable Taylor polynomials with generalized coefficients such that 
$(\tilde x, \tilde u^{(n)}) = \prn f (\tilde x)$. Therefore by definition
of the prolonged group action we have 
$$
L(\prn \Phi_{\tilde \eta}(\tilde x,  \tilde u^{(n)}))\det (J_{\tilde\eta}(\tilde x))
= L(\tilde x, \tilde u^{(n)})\,.
$$
Since this holds for any given  $(\tilde x, \tilde u^{(n)})\in ((\Omega_1\times \R^q)^{(n)})^\sim_c$, 
by the point value characterization of generalized functions we conclude that
\begin{equation} \label{eq:ic-5}
L(\prn \Phi_{\eta}(x,u^{(n)}))\det (J_{\eta}(x))
= L(x,u^{(n)}) \quad \mbox { in } \cG((-\eta_1,\eta_1)\times (\Omega_1\times\R^q)^{(n)})
\end{equation}
Differentiating this with respect to $\eta$ yields (by Lemma \ref{l:det J}):
\beq \label{eq:ic-4}
(\prn \vf(L) + L\Div \xi)\cdot \det (J_\eta(x))=0
\quad \mbox{in } \cG((-\eta_1,\eta_1)\times(\Omega_1\times\R^q)^{(n)}).
\eeq
If we set $\eta=0$, we obtain (\ref{eq:ic-1}) on $\Omega_1$. Since
$\Omega_1$ as above was arbitrary, the claim follows on all of $\Omega$.

Conversely, if (\ref{eq:ic-1}) holds in
$\cG((\Omega\times\R^q)^{(n)})$ for the infinitesimal generator $\xi$
of $\Phi$ then (\ref{eq:ic-4}) is automatically satisfied for each pair
$(\Omega_1,\eta_1)$ as above. Thus we may integrate (\ref{eq:ic-4}) from $0$ to 
$\eta$ to obtain (\ref{eq:ic-5}). Evaluating (\ref{eq:ic-5}) at the graph of
the generalized function $f\in \cG(\Omega_1)^q$ (cf.\ \cite[Definition 4.3.9]{book}), 
again by the point value characterization
in $\cG$ we obtain (\ref{eq:ic-3}), and, therefore, (\ref{eq:ic-2}), from which we 
conclude that $\Phi$ is a generalized variational symmetry group of $\cL$.
\ep

 As a preparation for N\"{o}ther's theorem,
 we briefly recall the prolongation formula in $\cG$ from \cite[Theorem 4.3.17]{book} 
 (based on \cite[Theorem 2.36]{Olv}): if $\vf = \sum_{i=1}^{p}\xi^i(x)\frac{\pd}{\pd x^i} +
 \sum_{\alpha=1}^{q}\psi_\alpha(x,u) \frac{\pd}{\pd u^\alpha}$ is
 a $\cG-n-$complete generalized vector field with corresponding
 projectable group action $\Phi$ on $\R^p\times\R^q$ then
 $$
 \prn\vf=\vf+\sum_{\alpha=1}^{q}\sum_{J} \psi^J_\alpha (x,
 u^{(n)}) \pd_{u^\alpha_J}, \,\,\, J=(j_1,\dots,j_k), 1\leq j_k
 \leq p, 1\leq k\leq n,
 $$
 where
 $$
 \psi^J_\alpha (x,\un)= D_J(\psi_\alpha - \sum_{i=1}^{p}
 \xi^iu_i^\alpha) + \sum_{i=1}^{p} \xi^i u^\alpha_{J,i} \in \cG_\tau((\R^p\times\R^q)^{(n)}).
 $$
 As in \cite[(2.48)]{Olv}, we may introduce the characteristics of $\vf$ as
 $$
 Q_\alpha(x,u^{(1)}):=\psi_\alpha(x,u) - \sum_{i=1}^{p}
 \xi^iu_i^\alpha \in \cG_\tau((\R^p\times\R^q)^{(1)}) \quad 1\leq\alpha\leq q,
 $$
 and $Q(x,u^{(1)}):=(Q_1,\dots,Q_q)(x,u^{(1)})$. Then
 $\psi_\alpha^J= D_JQ_\alpha + \sum_{i=1}^{p} \xi^iu^\alpha_{J,i}$
 and
 \beas
 \prn\vf &=& \sum_{\alpha=1}^{q}\sum_J D_JQ_\alpha \frac{\pd}{\pd
 u_J^\alpha}+ \sum_{i=1}^{p}\xi^i (\frac{\pd}{\pd x^i} +
 \sum_{\alpha=1}^{q}\sum_J u^\alpha_{J,i}\frac{\pd}{\pd u_J^\alpha})\\
 &=& \prn\vf_Q + \sum_{i=1}^{p} \xi^i D_i,
 \eeas
 with $D_i$ the $i$-th total derivative, and
 $$
 \vf_Q=\sum_{\alpha=1}^{q} Q_\alpha(x,u^{(1)}) \frac{\pd}{\pd
 u^\alpha},
 \quad \quad
 \prn\vf_Q = \sum_{\alpha=1}^{q}\sum_J D_J Q_\alpha
 \frac{\pd}{\pd u^\alpha_J}.
 $$

 N\"{o}ther's theorem provides a relation between variational
 symmetry groups and conservation laws for the corresponding
 system of Euler-Lagrange equations. Recall
 that a conservation law for a system of differential equations
 $\Delta(x,\un)=0$ is a divergence expression $\Div P=0$ which
 vanishes for all solutions $u=f(x)$ of $\Delta$ ($P$ is a
 $p$-tuple of generalized functions
 $P_i(x,\un)\in\cG((\R^{p+q})^{(n)})$, $1\leq i\leq p$).

 \bthm (N\"{o}ther's theorem) \label{th:Noether}
 Let $\Phi$ be a projectable, strictly
 $\cG$-complete variational symmetry group of the variational
 problem (\ref{eq:vp pq-dim}), and let
 $$
 \vf = \sum_{i=1}^{p}\xi^i(x)\frac{\pd}{\pd x^i} +
 \sum_{\alpha=1}^{q}\psi_\alpha(x,u) \frac{\pd}{\pd u^\alpha}
 $$
 be the infinitesimal generator of $\Phi$ with characteristics $Q_\alpha$ as above.
 Then $Q\cdot E(L)$ is a conservation law for the Euler-Lagrange
 equations $E(L)=0$, i.e., there exists a $p$-tuple $P(x, u^{(n)})
 =(P_1,\dots,P_p)(x, u^{(n)})$ such that 
 $$
 \Div P=Q\cdot E(L) = \sum_{\nu=1}^q Q_\nu E_\nu(L) \quad \mbox{in } \cG_\tau((\bar \Omega\times\R^q)^{(n)})\,.
 $$
 \ethm
 \pr
 The proof of the classical case (\cite[Theorem 4.29]{Olv}) carries over verbatim to our setting.
 \ep

For the practically most relevant case of first order variational problems, the explicit form of
$P$ is (cf.\ \cite[Cor.\ 4.30]{Olv}):
 $$
 P_i 
 = \sum_{\alpha=1}^q \psi_\alpha \frac{\pd L}{\pd u^\alpha_i}
 +\xi^i L -
 \sum_{\alpha=1}^q\sum_{j=1}^p \xi^j u_j^\alpha \frac{\pd L}{\pd
 u^\alpha_i}\,.
 $$

\section{Applications}
\label{sec:Appl}
\subsection{Geodesics as solutions of a variational problem}
\label{subsec:geodesics}
 Let $M$ be a smooth $q$-dimensional manifold equipped with a generalized Riemannian
 metric $g\in \cG^0_2(M)$ (cf.\ \cite{gprg,book}). We are looking for
 curves $c\in \cG[[0,T],M]$ which minimize the length functional
 $$
 L(c):=\int_0^T |\dot{c}(t)|\, dt = \int_0^T \sqrt{g(\dot c(t),\dot c(t))}\, dt
 $$
 As in the smooth setting (cf., e.g., \cite{Jost}), since $L(c) \le \sqrt{2T}\sqrt{E(c)}$ 
 it suffices to find minimizers of the energy functional
 \beq \label{eq:energy}
 E(c):=\frac{1}{2}\int_0^T |\dot{c}(t)|^2\, dt = \frac{1}{2}\int_0^T g(\dot c(t),\dot c(t))\, dt.
 \eeq
 In terms of a local coordinate system $(x^1,\dots,x^n)$, the Euler-Lagrange equations
 for $E$ read
 \begin{equation} \label{geoequ}
 \frac{d^2 c_\eps^k}{dt^2} + \sum_{i,j} \Gamma^k_{\eps ij} \frac{d c_\eps^i}{dt} \frac{d c_\eps^j}{dt} = 0
 \qquad k=1,\dots,q
 \end{equation}
 with $\Gamma^k_{ij}$ the Christoffel symbols of $g$,
 $$
 \Gamma^k_{ij} = \frac12 \sum_m g^{km}\left(\frac{\partial g_{jm}}{\partial x^i} + 
  \frac{\partial g_{im}}{\partial x^j} 
 - \frac{\partial g_{ij}}{\partial x^m} \right)\,,
 $$
 $g_{ij} =  g(\frac{\partial}{\partial x^i},\frac{\partial}{\partial x^j})$, and $(g^{ij}) = (g_{ij})^{-1}$.
 These are the usual geodesic equations for the metric $g$, whose derivation from (\ref{eq:energy}) in the
 Colombeau setting, based on Lemma \ref{l:fundamental pq-dim}, is completely analogous to that 
 in the smooth case (\cite{Jost}). (\ref{geoequ}) is a system of ordinary
 differential equations for the components of $c$ in $\cG([0,T])$.

 Analogously, for pseudo-Riemannian generalized metrics, (\ref{geoequ}) is derived as the Euler-Lagrange equation
 for $E$ by fixing the causal character (cf.\ \cite{EbM-AF}) of $c$.
 
 As a concrete example let us consider the case of impulsive pp-waves ({\em plane fronted gravitational waves 
 with parallel rays}) (\cite{herbertgeo, geo, geo2}). Here, the metric is given in the form
 $$
 g = f(x,y) \delta(u) du^2 - dudv +dx^2 + dy^2
 $$ 
 with $f$ smooth and $\delta$ the Dirac measure.
 Using a generic regularization of $\delta$ to embed $g$ into the Colombeau algebra, 
 one can show that the corresponding
 system of geodesic equations is uniquely solvable in $\cG[\R]$. Moreover, this unique solution is associated
 to the physically expected solution (piecewise straight geodesics which are broken and refracted at the 
 shock hypersurface $u=0$). We refer to the above papers or also to \cite[Sec.\ 5.3]{book} for a detailed
 exposition.

\subsection{Examples from mechanics}
\label{subsec:mechanics}
Denote by $x(t)$ the position of a particle of mass 1 in $\R^n$ that moves in
a potential $V(x)$. According to Hamilton's least action principle,
the particle follows the path that minimizes the action functional
\begin{equation}\label{eq:action}
   \cL(x) = \int_a^b L(x(t),\dot{x}(t))\,dt
\end{equation}
among all trajectories $y(t)$ such that $y(a) = x(a)$, $y(b) = x(b)$
for whatever points of time $a < b$. Here the Lagrangian is of the
form
\begin{equation}\label{eq:T-U}
   L(t,x,\dot{x}) = \frac{1}{2} |\dot{x}|^2 - V(x)\;
\end{equation}
and given by the difference of the kinetic and the potential energy.
It is clear that time translations form a variational symmetry for
the functional \eqref{eq:action}. It follows from
N\"other's theorem that the total energy
\[
   P = \frac{1}{2} |\dot{x}|^2 + V(x)
\]
is conserved. Particle systems are of interest here, because our
approach offers the possibility to allow singular potentials, for
example, delta function potentials. Singular potentials may be modelled
by elements $V\in\cG_\tau(\R^n)$. If $x\in (\cG(\R))^n$ minimizes
the action functional \eqref{eq:action} on every interval $[a,b]$,
then it follows from Theorem \ref{eulerlagrange} that $x$ satisfies the
Euler-Lagrange-equations
\begin{equation}  \label{eq:EulerLagrange}
   \ddot{x}(t) + \nabla V(x) = 0\,.
\end{equation}
It has been shown in \cite{book} that \eqref{eq:EulerLagrange}, completed by initial data
$x(t_0) = \widetilde{x}_0$, $\dot{x}(t_0) = \widetilde{y}_0$ has a
solution $x\in (\cG(\R))^n$, provided $V$ is either of $L^\infty$-type or nonnegative.
In addition, total energy is conserved along the generalized trajectories.

Uniqueness of the solution to (\ref{eq:EulerLagrange}) in $(\cG(\R))^n$
can be obtained when $\nabla^2 V$ is of $L^\infty$-log-type.
As in the classical case, the second variation of \eqref{eq:action} is
rarely positive so that uniqueness cannot be inferred using the
variational principle alone, in general.

We now turn to a more explicit example -- a classical particle in a delta function
potential. Formally, it is described by the Lagrangian
$L(t,x,\dot{x}) = \frac{1}{2}|\dot{x}|^2 - \delta(x)$.
The delta function potential is considered as acting as a barrier so that the particle
trajectories exhibit pure reflection at the origin.
To make things more precise, let us call
a generalized function $D \in \cG_\tau(\R)$ a {\em strict delta function}
if it has a representative $(\rho_\eps)_\eps \in \cE_\tau(\R)$ such that
\begin{itemize}
\item[(i)] $\supp\, \rho_\eps \to \{0\}$ for $\eps\to 0$.
\item[(ii)] $\lim_{\eps\to 0} \int \rho_\eps(x) dx =1$.
\item[(iii)] $\int |\rho_\eps(x)|\, dx$ is bounded uniformly in
$\eps$.
\end{itemize}
$D$ is called a {\em model delta function}
if it has a representative $(\varphi_\eps)_\eps$ of the form
$\varphi_\eps(x) = \frac{1}{\eps} \varphi(\frac{x}{\eps})$,
where $\varphi\in {\cal D}(\R)$, $\int\varphi(x)\, dx=1$.
%
\bex
(Particle in a delta function potential)
The Euler-Lagrange equations are
\begin{eqnarray}
\ddot{x}(t) + D'(x(t)) &=& 0 \nonumber\\
                   x(0)&=&\widetilde{x}_0  \label{eq:deltode} \\
             \dot{x}(0)&=&\widetilde{y}_0\nonumber
\end{eqnarray}
where $D \in \cG_\tau(\R)$ is a strict or model delta function
and $\widetilde{x}_0,\widetilde{y}_0\in \tilde{\R}$. We have the following result:
\begin{itemize}
\item[(i)] Problem (\ref{eq:deltode}) has a solution $x\in \cG(\R)$.
\item[(ii)] If $D''$ is of $L^\infty$-log-type, the solution is unique in $\cG(\R)$.
\item[(iii)] If $D$ is a model delta function and
$\widetilde{x}_0=x_0,\widetilde{y}_0=y_0$ are elements of $\R$ with
$x_0\not= 0$, then the solution $x \in \cG(\R)$ of (\ref{eq:deltode}) is associated
with the function $t\mapsto {\rm sign}(x_0)|x_0+t y_0|$.
\end{itemize}
\eex
The pure reflection picture does not necessarily apply when strict delta functions
are employed. In general, it is possible that certain trajectories terminate at the origin,
so that the particle may be trapped there after finite time. More details and proofs
of these results can be found in \cite{book}.

\bex
(Planar motion in a central force field)
Introducing polar coordinates $x = r\cos\varphi$, $y = r\sin\varphi$ in the plane,
the motion of a particle of mass $m$ in the time interval $[a,b]$ is a minimizer of the action functional
\[
   \cL(r,\varphi) = \int_a^b L(r(t),\dot{r}(t),\varphi(t),\dot{\varphi}(t))\,dt
\]   
with
\begin{equation}\label{eq:centralfield}
   L(r,\dot{r},\varphi,\dot{\varphi}) = \frac{1}{2} m\big(\dot{r}^2 + r^2\dot{\varphi}^2\big) - V(r).
\end{equation}
Here the potential energy $V(r)$ is assumed to be an element of $\cG_\tau(\R)$.
Such a situation could arise from distributional potentials as in the previous example,
or from regularizing classical, non-smooth potentials like the gravitational
potential $-1/r$. It is clear that the infinitesimal generator of the rotation
\[
   \vf = \frac{\p}{\p\dot{\varphi}}
\]
determines a variational symmetry. According to N\"other's theorem (Theorem \ref{th:Noether})
the corresponding conserved quantity is angular momentum 
\[
  P = \frac{\p L}{\p\dot{\varphi}} = mr^2\dot{\varphi},
\]
and this holds for generalized solutions $r,\varphi \in \cG[a,b]$ as well.
\eex

\subsection{Examples from elastostatics}
\label{subsec:elastostatics}
We begin with an example from classical elastostatics. To keep matters
simple, we do not treat the full three-dimensional problem but rather focus on
a model problem that describes, e.\;g., vertical displacements of membranes ($n=2$)
or strings ($n=1$).
\bex
(Standard linear elasticity)
Let $\Omega\subset \R^n$ be a bounded,
open and connected set with smooth boundary $\partial\Omega$, the reference configuration
of an elastic body. We denote by $u: \bar{\Omega} \to \R$ the displacement.
For simplicity, we address the Dirichlet problem $u|_{\partial\Omega} = u_0$
only. The static solution minimizes the functional
\[
   \cL(u) = \frac{1}{2}\int_\Omega \big(\alpha(x)|\nabla u(x)|^2 + \beta(x)u(x)^2\big)\,dx
                     - \int_\Omega \gamma(x)u(x)\,dx
\]
\eex
under all displacements that satisfy the boundary condition. Here $\alpha(x)$ depends on the
elastic properties of the material (modulus of elasticity and Poisson's ratio),
$\beta(x)$ is a potential that could arise by an elastic
constraint, for example, and $\gamma(x)$ is the body force. In our setting, we may take
$\alpha$, $\beta$ and $\gamma$ as generalized functions. In this way, singular potentials and degenerate
coefficients $\alpha(x)$ can be modelled. A delta potential $\beta(x) = \delta(x-x_0)$ would arise
if a spring is attached to a membrane at point $x_0\in\Omega$. If the coefficient $\alpha(x)$
vanishes in a subregion of $\Omega$, the classical problem looses wellposedness, but the
generalized problem might still be wellposed if $\alpha(x)$ is modelled by a generalized function
associated with zero, but still invertible.

As in Example~\ref{ex:quadraticforms}, we make the following assumptions: 
$\alpha, \beta, \gamma \in \cG(\bar{\Omega})$, $\alpha \geq 0$ and invertible
on $\bar{\Omega}$, $\beta \geq 0$, $u_0 \in \cG(\p\Omega)$. The admissible set and the
admissible variations, respectively, are given by
\[
   \cA = \{u \in \cG(\bar{\Omega}): u|_{\p\Omega} = u_0\},\qquad
        \cA_0 = \{u \in \cG(\bar{\Omega}): u|_{\p\Omega} = 0\}.
\]
Observe that $\cL$ is the quadratic functional
\[
   \cL(u) = \frac{1}{2}a(u,u) - f(u)
\]
corresponding to the the bilinear form
\[
   a(u,v) = \int_\Omega \big(\alpha(x)\nabla u(x)\cdot\nabla v(x) + \beta(x)u(x)v(x)\big)\,dx
\]
and the linear form
\[
   f(u) = \int_\Omega \gamma(x)u(x)\,dx.
\]
The Euler-Lagrange equation is
\begin{equation}\label{eq:ELelasticity}
    -\DIV \big(\alpha(x)\nabla u(x)\big) + \beta(x) u(x) = \gamma(x), \qquad u|_{\p\Omega} = u_0.
\end{equation}
The first and second variation was calculated in Example~\ref{ex:quadraticforms} and
it was shown that $\cL(u)$ admits a unique minimizer.
It is clear that $u \in \cA$ solves the Euler-Lagrange equation
\eqref{eq:ELelasticity} in $\cG(\Omega)$ if and only if it solves the
variational equation (\ref{VE}).
In this way, we have proved that solutions to the Euler-Lagrange equation are unique.

In case the coefficients $\alpha, \beta, \gamma$ are classical smooth functions,
a proof of existence can be obtained as follows: Fixing a representative $u_{0\eps}$
of $u_0$, equation \eqref{eq:ELelasticity} has a unique classical smooth
solution $u_\eps$, $\eps\in (0,1]$, cf. \cite[Chap.\;6]{Gilbarg}. 
The closed graph theorem yields continuous dependence on the boundary data in the
$\cC^\infty$-seminorms. This shows that $(u_\eps)_\eps$ is
moderate and thus defines a solution $u\in \cG(\bar{\Omega})$ to \eqref{eq:ELelasticity}.
An existence result with generalized coefficients $\alpha\in \cG(\bar{\Omega})$,
$\beta = \gamma = 0$ has been obtained in \cite{Pilipovic-Scarpalezos}.
Alternatively, the solution can be constructed by direct methods of the calculus of
variations, as have been worked out in the Colombeau setting in \cite{GarettoVernaeve}.

\bex
(Euler-Bernoulli beam with discontinuous coefficient)
The equilibrium equation for the deflection $w(x)$, $x\in [0,\ell]$ of an Euler-Bernoulli beam
under constant axial force $\beta$ and distributed load $\gamma(x)$ is
\begin{equation}\label{eq:E-B}
   \frac{d^2}{dx^2}\Big(\alpha(x)\frac{d^2w(x)}{dx^2}\Big) + \beta \frac{d^2w(x)}{dx^2} - \gamma(x) = 0,\quad x\in [0,\ell].
\end{equation}
\eex
Here $\ell$ is the length of the beam. The coefficient $\alpha$ is given by the flexural stiffness
$EI$ ($E$ the modulus of elasticity, $I$ the area moment of inertia) 
and may depend on the position $x$ along the beam, in general.
Appropriate boundary conditions for a beam freely supported at both ends are
\begin{equation}\label{eq:beamBC}
   w(0) = 0,\quad w(\ell) = 0,\quad \frac{d^2}{dx^2}w(0) = 0,\quad \frac{d^2}{dx^2}w(\ell) = 0.
\end{equation}
We are particularly interested in the following cases:
\begin{itemize}
\item[(a)] the coefficient $\alpha$ has jump discontinuities (non-smooth case);
\item[(b)] the coefficient $\alpha$ vanishes in certain points (degenerate case);
\item[(c)] the coefficient $\alpha$ is given by a path of positive noise.
\end{itemize}
Case (a) is a simple model of abrupt change of material properties (see e.g. \cite{HoermannOparnica, YavariSarkani});
case (b) arises for example, if the beam contains a joint in its interior. Classically, these
problems are split into adjacent subproblems with transition conditions. A global formulation
in the classical setting -- particularly in case (b) -- is impossible, because the classical
ellipticity condition would be violated. Case (c) is of relevance in stochastic structural
mechanics. Since $\alpha$ has to be nonnegative, a lognormally distributed noise process would be a suitable
model for random disturbances. Such a process does not exist in the classical theory of stochastic processes,
but it does exist in the sense of a Wick exponential \cite{OksendalSPDEs} or as a Colombeau
process \cite{MORajter}. In all three cases the Colombeau setting lends itself as a framework that 
provides a rigorous global solution concept.

We observe that (\ref{eq:E-B}) is the Euler-Lagrange equation corresponding to the quadratic functional
\[
   \cL(u) = \frac{1}{2}\int_0^\ell \big(\alpha(x)|u''(x)|^2 -  \beta|u'(x)|^2\big)\,dx
                     - \int_0^\ell \gamma(x)u(x)\,dx.
\]
The second variation is given by
\begin{equation}\label{eq:E-B2Var}
  \delta^2\cL(u,v) = \frac{1}{2}\int_0^\ell \big(\alpha(x)|v''(x)|^2 -  \beta|v'(x)|^2\big)\,dx
\end{equation}
We assume that there is an invertible element $\alpha_0 \geq 0$ in $\tilde{\R}$ such that
$\alpha \geq \alpha_0$ in $\cG[0,\ell]$. Clearly, the second variation (\ref{eq:E-B2Var}) is
positive definite if the axial force $\beta$ is negative. We are going to show that
(\ref{eq:E-B2Var}) is positive definite even for small, positive axial forces $\beta$, 
below the first eigenvalue
of equation (\ref{eq:E-B}). To see this, let $w\in \cC^\infty[0,\ell]$, $w(0) = w(\ell) = 0$. Then
\[
  w(x) = -\frac{x}{\ell}\int_0^\ell (\ell - y)w''(y)\,dy + \int_0^x(x-y)w''(y)\,dy.
\] From there we derive the inequality
\[
   \|w\|_{{\rm L}^2(0,\ell)}^2 \leq C\|w''\|_{{\rm L}^2(0,\ell)}^2
\]
for some constant $C > 0$ and all smooth $w$ as above. Further,
\begin{eqnarray*}
   - \int_0^\ell|w'(x)|^2\,dx &=& \int_0^\ell w''(x)w(x)\,dx\\
      &\geq& -\frac{1}{2} \int_0^\ell |w''(x)|^2\,dx -\frac{1}{2} \int_0^\ell |w(x)|^2\,dx\\
      &\geq& -\frac{1}{2}(1+C) \int_0^\ell |w''(x)|^2\,dx.
\end{eqnarray*} From this inequality it follows that (\ref{eq:E-B2Var}) is positive definite if
$\alpha_0 - \frac{1}{2}(1+C)\beta \geq 0$ in $\cG[0,\ell]$. If this expression 
is invertible, then $\cL$ admits a unique minimizer (Prop.~\ref{prop:cv-sufficient cond}),
which is obtained as the solution to the Euler-Lagrange equation (\ref{eq:E-B}).

As a concrete application, we are going to work out a global formulation of an
Euler-Bernoulli beam with a joint, say at half-length $x = \ell/2$. At a joint,
the flexural stiffness , i.e. the coefficient $\alpha$, vanishes. In the
Colombeau setting, this can be modelled as follows. Take a smooth unimodal
function $\psi$ vanishing for $|x| \geq 1$ and equal to $1$ for $|x|\leq 1/2$,
$0 \leq \psi(x) \leq 1$ and let $h\approx 0$ define a nonnegative invertible element of $\tilde{\R}$. Put
\[
   \alpha_\eps(x) = \alpha \Big(1 - \big(1-h_\eps\big)\psi\Big(\frac{x-\ell/2}{\eps}\Big)\Big)
\]
where $\alpha\in \R$ is some positive constant. Clearly, $(\alpha_\eps)_\eps$ determines a nonnegative 
invertible element of $\cG([0,\ell]$ which is infinitesimally small near $x = \ell/2$ and otherwise equal
to $\alpha$.
We consider the Euler-Bernoulli beam without axial force
\begin{equation}\label{eq:E-B0}
   \frac{d^2}{dx^2}\Big(\alpha_\eps(x)\frac{d^2u_\eps}{dx^2}(x)\Big) - \gamma(x) = 0,\quad x\in [0,\ell]
\end{equation}
where $\gamma \in \cC^\infty[0,\ell]$ denotes the distributed load.
A representative $u_\eps$ of the solution in $\cG[0,\ell]$ is readily calculated as
\[
   u_\eps(x) = \int_0^x\int_0^y \frac{M(z)}{\alpha_\eps(z)}\,dzdy 
     - \frac{x}{\ell} \int_0^\ell\int_0^y \frac{M(z)}{\alpha_\eps(z)}\,dzdy 
\]
where $M(x)$, minus the bending moment, is the solution to $M''(x) = \gamma(x)$, $M(0) = M(\ell) = 0$.
Now,
\[
   \int_0^y \frac{M(z)}{\alpha_\eps(z)}\,dz = \left\{
   \begin{array}{ll}
   \int_0^y \frac{M(z)}{\alpha}\,dz, & 0 \leq y \leq \ell/2 - \eps,\\[6pt]
   \int_0^{\ell/2-\eps} \frac{M(z)}{\alpha}\,dz + C_\eps + \int_{\ell/2+\eps}^y \frac{M(z)}{\alpha}\,dz,
			& \ell/2+\eps\leq y \leq \ell
   \end{array}\right.
\]
with
\[
   C_\eps = \int_{\ell/2-\eps}^{\ell/2 + \eps} \frac{M(z)}{\alpha_\eps(z)}\,dz.
\]
If we choose $h_\eps$ in such a way that 
\begin{equation}\label{eq:D}
   \int_{\ell/2-\eps}^{\ell/2 + \eps} \frac{dz}{1 - (1-h_\eps)\psi\big(\frac{z-\ell/2}{\eps}\big)}
      = \int_{-1}^1 \frac{\eps\, dy}{1 - (1-h_\eps)\psi(y)} \to D
\end{equation}
in $\R$ as $\eps \to 0$, we will have that
\[
   \lim_{\eps \to 0}\int_0^y \frac{M(z)}{\alpha_\eps(z)}\,dz = \left\{
   \begin{array}{ll}
   \int_0^y \frac{M(z)}{\alpha}\,dz, & 0 \leq y \leq \ell/2,\\[6pt]
   \int_0^{\ell/2} \frac{M(z)}{\alpha}\,dz + \alpha DM(\frac{\ell}{2}) + \int_{\ell/2}^y \frac{M(z)}{\alpha}\,dz,
			& \ell/2 \leq y \leq \ell
   \end{array}\right.
\]
and so $u_\eps(x)$ will converge to a continuous, piecewise smooth function $u(x)$. 
This limit -- the associated distribution --  describes the displacement curve of the beam. 
Larger values of $D$ correspond to larger loss of stiffness at $x = \ell/2$ and thus to
larger displacements of center of the beam. Failure of the beam can be modelled by letting
$D = \infty$ in (\ref{eq:D}).

We remark that the case of a discontinuous flexural stiffness and nonzero axial force
has been treated in \cite{HoermannOparnica}. Finally, in the stochastic case, positive noise on $[0,\ell]$ 
can be defined in the Colombeau sense as
a lognormally distributed process with mean value 1 and variance given by 
$\exp(\|\rho_\eps\|_{{\rm L}^2(\R)}^2)$ $-$ $1$,
where $\rho_\eps$ is a strict delta function as in Subsection~\ref{subsec:mechanics}.
The paths of this generalized random process are nonnegative elements of the 
Colombeau algebra $\cG[0,\ell]$ and thus may serve as describing highly random behavior
of the flexural stiffness $\alpha$. We refer to \cite{MORajter} for details on and
explanations of the positive noise model.

\bex
(Rods with generalized stress-strain relation; hard rods)
Let $u(x)$ be the displacement of a rod of length $\ell$ of cross-sectional
area one, subject to the body force (density) $f(x)$, $0\leq x \leq \ell$. The balance law is
\[
   \bsig'(x) + f(x) = 0
\]
where $\bsig$ denotes stress. Letting $\beps = \frac{\p u}{dx}$ the strain, assume a
stress-strain relation (constitutive law) of the form $\bsig = g(\beps)$. If the rod is clamped
at one end and free at the other, the displacement is the solution to the
problem
\begin{equation}\label{eq:rod}
\frac{d}{dx}\, g\Big(\frac{d}{dx}u(x)\Big) + f(x) = 0, \quad u(0) = 0, \ g\Big(\frac{d}{dx}u(\ell)\Big) = 0.
\end{equation}
If the constitutive law has a potential, $g(y) = -G'(y)$, equation (\ref{eq:rod}) is the
Euler-Lagrange equation of the functional
\begin{equation}\label{eq:rodL}
  \cL(u) = \int_0^\ell \Big( G(u'(x)) + f(x)u(x)\Big) dx.
\end{equation}
We allow generalized potentials $G \in \cG_\tau[0,\ell]$ and forces $f \in \cG[0,\ell]$.
The admissible set and admissible variations are
\[
   \cA = \cA_0 = \{u\in \cG[0,\ell]: u(0) = 0\}.
\]   
Assume that $u\in\cA$ minimizes the functional (\ref{eq:rodL}). The first variation is
\begin{eqnarray*}
   \delta\cL(u;v) & = & \int_0^\ell \Big(-g(u'(x))v'(x) + f(x)v(x)\Big) dx \\
      & = & \int_0^\ell \Big(g(u'(x))'v(x) + f(x)v(x)\Big) dx + g(u'(\ell))v(\ell).
\end{eqnarray*}
If $\delta\cL(u;v) = 0$ for all $v\in \cA_0$, we have in particular that
\[
   g(u'(x))' + f(x) = 0
\]
in $\cG(0,\ell)$, by using Lemma~\ref{l:fundamental pq-dim}. In order to show that 
$g(u'(\ell)) = 0$ we need an adaptation of the proof of this lemma. Thus assume that
$g(u'(\ell)) \neq 0$ in $\tilde{\R}$. We can find $(r_\eps)_\eps \in \R_M$ and a subsequence
$\eps_k \to 0$ such that $g_{\eps_k}(u_{\eps_k'}(\ell))\,r_{\eps_k} = 1$ for all $k\in\N$.
On the other hand, there is $N \geq 0$ such that
\[
  \sup_{0 \leq x \leq \ell}|g_\eps(u_\eps'(x))' + f_\eps(x)| = O(\eps^{-N})\ {\rm as}\ \eps\to 0.
\]
Similar to the argument in Lemma~\ref{l:fundamental pq-dim}, we can find an element
$v\in \cG[0,\ell]$ with support in $(0,\ell]$ such that $v_{\eps_k}(\ell) = r_{\eps_k}$ and
\[
  \Big|\int_0^\ell \big(g_\eps(u_\eps'(x)) + f_\eps(x)\big)v_\eps(x) dx\Big| \leq \frac12
\]
for all $\eps$, contradicting the hypothesis that $\delta\cL(u;v) = 0$ for all $v\in \cA_0$.
We arrive at the conclusion that a minimizer of (\ref{eq:rodL}) in $\cA$ is a solution to
the Euler-Lagrange equation (\ref{eq:rod}).

In linear elasticity, the stress-strain relation is $\bsig = E\beps$ with the modulus of elasticity
$E$. The larger $E$, the harder is the rod. In the limiting case $E\to\infty$, the rod becomes
inextensible. We can model such a rod in the Colombeau framework by introducing
the potential $G\in  \cG_\tau(\R)$ by $G_\eps(y) = -\frac{1}{\eps}y$.
The Euler-Lagrange equation becomes 
\[
\frac{1}{\eps}u_\eps''(x) + f_\eps(x) = 0, \quad u_\eps(0) = 0, \ u_\eps'(\ell) = 0
\]
and has the solution
\[
   u_\eps(x) = \eps \int_0^x\int_y^\ell f_\eps(z)\,dzdy\,.
\]
For bounded body force, the solution is associated with zero, thus no extension of the
rod takes place, as anticipated. We refer to \cite{Moreau} for the corresponding
model in convex analysis.
\eex

\subsection{Nonlinear wave equations with singular potential}
\label{subsec:singpot}
\bex
(Wave equation with delta potential)
In this subsection, we discuss the one-dimen\-sio\-nal wave
equation 
\begin{equation}\label{eq:wave}
   \frac{\p^2}{\p t^2} u(x,t) - \frac{\p^2}{\p x^2} u(x,t) + \frac{\p}{\p u} V\big(x,u(x,t)\big) = 0
\end{equation}
with potential $V(x,u)$. The higher dimensional case can be treated similarly.
Classical examples of the potential are $V(u) = \frac{m}{2}u^2 + \frac{k}{2}u^4$, $m,k > 0$,
leading to the cubic Klein-Gordon equation, or $V(u) = - \cos u$, leading to the
Sine-Gordon equation. However, much more singular potentials are in use. For example,
the equation
\begin{equation}\label{eq:deltawave}
   \frac{\p^2}{\p t^2} u(x,t) - \frac{\p^2}{\p x^2} u(x,t) + \delta(x-x_0)F(u(x,t)) = 0
\end{equation}
describes the vibrations of a string with a nonlinear spring with restoring force $F(u)$
attached at the point $x = x_0$, see \cite{Komech}. An even more singular potential,
\[
   \frac{\p^2}{\p t^2} u(x,t) - \frac{\p^2}{\p x^2} u(x,t) - \sum_{k=1}^n\delta(u(x,t) - u_k) = 0
\]
where $\{u_1, \ldots, u_n\}$ is a finite subset of $\R$, has been studied in \cite{Bensoussan}.
This equation arises from a piecewise constant approximation of a smooth potential -- in the quoted reference,
the Ginzburg-Landau potential. 
Note that in this example, the sought after solution has to be inserted in a sum of Dirac measures.
All these cases can be subsumed by letting the potential $V(x,u)$ belong to the Colombeau algebra 
$\tilcG_\tau(\R\times\R)$. The Lagrangian functional corresponding to (\ref{eq:wave}) is
\[
  \cL(u) = \int_a^b\int_c^d \Big(\frac12 \Big|\frac{\p}{\p t}u(x,t)\Big|^2
	- \frac12 \Big|\frac{\p}{\p x}u(x,t)\Big|^2 - V\big(x, u(x,t)\big)\Big)dxdt.
\]
The element $u \in \cG(\R^2)$ that makes $\cL(u)$ stationary on every rectangle $[a,b]\times[c,d]$
produces a solution of the Euler-Lagrange equation (\ref{eq:wave}).

Existence and uniqueness of solutions in $\cG(\R^2)$ to the Cauchy problem for equation (\ref{eq:wave})
have been demonstrated under various conditions on the potential $V$, see e.g. \cite{MONA, MONedelPili}.
Solutions to equation (\ref{eq:deltawave}) in the Colombeau algebra $\cG(\R^2)$ and their associated
distributions have been computed in \cite{MOwaves}.
\eex



\end{document}